\documentclass[preprint,3p,sort&compress,final,times]{elsarticle}
\usepackage{amsmath}
\usepackage{amssymb}
\usepackage{mathtools}
\usepackage{placeins}
\usepackage[usenames,dvipsnames]{color}
\usepackage{xcolor}
\usepackage{cases}
\usepackage{array}

\usepackage[usenames,dvipsnames]{color}

\newcommand{\firstRev}[1]{\color{black}{#1}\color{black}}
\newcommand{\secondRev}[1]{\color{black}{#1}\color{black}}

\newcommand{\bm}[1]{\boldsymbol{\mathsf{#1}}}            % a matrix operator
\newcommand{\ip}[2]{\langle #1, #2\rangle_{\Omega}}

\journal{Journal of Computational Physics}
\begin{document}
\begin{frontmatter}

\title{A Mixed Mimetic Spectral Element Model of the Rotating Shallow Water Equations on the Cubed Sphere}
\author[MON,CCS]{D.~Lee\corref{cor}}
\ead{davidr.lee@monash.edu}
\author[EUT]{A.~Palha}

\address[MON]{Department of Mechanical and Aerospace Engineering, Monash Univeristy, Melbourne 3800, Australia}
\address[CCS]{Computer, Computational and Statistical Sciences, Los Alamos National Laboratory, Los Alamos, NM 87545, USA}
\address[EUT]{Eindhoven University of Technology, Department of Mechanical Engineering, P.O. Box 513, 5600 MB Eindhoven, 
The Netherlands}

\cortext[cor]{Corresponding author. Tel. +61 452 262 804.}

\begin{abstract}
In a previous article [\emph{J. Comp. Phys.} $\mathbf{357}$ (2018) 282-304], the mixed mimetic 
spectral element method was used to solve the rotating shallow water equations in an idealized geometry. Here the 
method is extended to a smoothly varying, non-affine, cubed sphere geometry. The differential operators are encoded 
topologically via incidence matrices due to the use of spectral element edge functions to construct tensor product 
solution spaces in $H(\mathrm{rot})$, $H(\mathrm{div})$ and $L_2$. These incidence matrices commute with respect 
to the metric terms in order to ensure that the mimetic properties are preserved independent of the geometry. This 
ensures conservation of mass, vorticity and energy for the rotating shallow water equations using inexact 
quadrature on the cubed sphere. The spectral convergence of errors are similarly preserved on the cubed sphere, 
with the generalized Piola transformation used to construct the metric terms for the physical field quantities.
\end{abstract}

\begin{keyword}
Mimetic\sep
Spectral convergence\sep
Shallow water\sep
Cubed sphere
\end{keyword}

\end{frontmatter}

\section{Introduction}\label{sec::intro}

In recent years there has been much attention given to the use of mimetic or compatible finite element methods for the
modelling of geophysical flows. This work has been motivated by the desire to preserve conservation laws in order to 
mitigate against biases in the solution over long time integrations \cite{Thuburn08}. 
These mimetic methods are designed to preserve the divergence and circulation theorems in the discrete form, as well as 
the annihilation of the gradient by the curl and the curl by the divergence. When appropriate solution spaces are chosen 
for the divergent, vector and rotational moments, this allows for the conservation of first (mass, vorticity) and higher 
(energy and potential enstrophy) moments in the discrete form \cite{AL81,MC14,LPG18}. Various classes of element types
have been explored for this purpose, including Raviart-Thomas, Brezzi-Douglas-Marini and Brezzi-Douglas-Fortin-Marini 
elements \cite{CS12,MC14,NSC16,SC17}. Mimetic properties may also be recovered for standard collocated A-grid spectral 
elements \cite{TF10} and primal/dual finite volume formulations \cite{RTKS10}.

When implemented on non-affine geometries, the convergence of errors \firstRev{ may }
%has been observed to 
degrade for compatible finite element methods \cite{ABB15}, \firstRev{ due to the reduced order of the 
polynomials when scaled by non-constant metric terms}. 
Several methods have been shown to rehabilitate the optimal convergence of Raviart-Thomas
elements \firstRev{for the $L_2$ function space } \cite{BR08,BG09,SC17} by modifying how the metric terms are incorporated 
into the differential operators, however it is unclear if and how these methods are applicable to other families of compatible 
finite element methods.

In the present article we extend previous work on the use of mixed mimetic spectral elements for geophysical flows
\cite{LPG18}, hereafter LPG18, to a non-affine cubed sphere geometry. The method uses the spectral element edge functions
\cite{Gerritsma11}, which are specified to satisfy the Kronecker delta property with respect to their integrals between 
nodes, so as to exactly satisfy the fundamental theorem of calculus with respect to the standard nodal spectral element 
basis functions. Combinations of standard nodal and edge functions are then used to construct tensor product solution 
spaces in higher dimensions for which the differential operators may be defined in a purely topological manner via the 
use of incidence matrices \cite{KPG11,KG13,PRHKG14}. These incidence matrices allow for the preservation of the divergence 
and circulation theorems, as well as the annihilation of the gradient by the curl and the curl by the divergence in the 
discrete form. The incidence matrices also commute with the metric transformations between computational and physical space, 
such that both the mimetic properties and the spectral convergence of errors are preserved 
\secondRev{ on smoothly varying, non-affine geometries } \cite{KPG11,PRHKG14}. \secondRev{ Indeed, for the spectral mimetic
least squares method, optimal convergence has also been demonstrated for irregular meshes that do not vary smoothly or 
converge to an affine geometry } \cite{BG14}.

In LPG18 the conservation and convergence properties of the mixed mimetic spectral element method for rotating shallow water 
flows were demonstrated both theoretically through formal proofs in the discrete form, and experimentally, through numerical 
experiments on idealized doubly periodic geometries. Here we extend these results to a non-affine cubed sphere geometry via 
the use of the generalized Piola transformation \cite{RKL09,RHCM13,NSC16}. This demonstrates that both the conservation laws 
derived from the mimetic properties, and the spectral convergence of errors, are preserved 
\secondRev{ for the smoothly varying, non-affine mesh of the cubed sphere }
without the need to rehabilitate the method through the modification of the discrete differential operators.

The remainder of this article proceeds as follows. In Section \ref{sec::sem} the formulation of the mixed mimetic spectral 
element method will be briefly discussed. Section \ref{sec::met} will discuss the formulation of the metric terms and their 
commuting properties with respect to the differential operators. The solution of the rotating shallow water equations on the
cubed sphere using mixed mimetic spectral elements will be discussed in Section \ref{sec::sw}. Section \ref{sec::res} will 
present results from some standard test cases demonstrating the preservation of optimal spectral convergence and conservation 
laws on the cubed sphere, and finally Section \ref{sec::conc} will discuss the conclusions of this work and some future 
directions we intend to pursue with this research.

\section{Mixed mimetic spectral elements}\label{sec::sem}

In this section we introduce the construction of the mixed mimetic spectral element method. For a more detailed discussion 
see LPG18, as well as previous work \cite{Gerritsma11,KPG11,KG13,PRHKG14} and references therein.

\subsection{One dimensional nodal and histopolant polynomials}

The mixed mimetic spectral element method is built off two types of one-dimensional polynomials: one associated with nodal 
interpolation, and the other with integral interpolation (histopolation) \cite{robidoux-polynomial,Gerritsma11}. Subsequently, 
these two types of polynomials will be combined to generate the family of two-dimensional polynomial basis functions used to 
discretize the system.

Consider the canonical interval $I=[-1,1]\subset\mathbb{R}$ and the Legendre polynomials, $L_{p}(\xi)$ of degree $p$ with 
$\xi\in I$. The $p+1$ roots, $\xi_{i}$, of the polynomial $(1-\xi^{2})\frac{\mathrm{d}L_{p}}{\mathrm{d}\xi}$ are called 
Gauss-Lobatto-Legendre (GLL) nodes and satisfy $-1 = \xi_{0} < \xi_{1} < \dots < \xi_{p-1} < \xi_{p} = 1$. Let $l^{p}_{i}(\xi)$ 
be the Lagrange polynomial of degree $p$ through the GLL nodes, such that
            \begin{equation}
                l^{p}_{i}(\xi_{j}) :=
                \begin{cases}
                    1 & \mbox{if } i = j \\
                    & \\
                    0 & \mbox{if } i \neq j
                \end{cases}\,, \quad i,j = 0,\dots,p\,. \label{eq::nodal_basis_polynomials}
            \end{equation}
            The explicit form of these Lagrange polynomials is given by
            \begin{equation}
                l^{p}_{i}(\xi) = \prod_{\substack{k=0\\k\neq i}}^{p}\frac{\xi-\xi_{k}}{\xi_{i}-\xi_{k}}\,. \label{eq::lagrange_interpolants}
            \end{equation}

Let $q_h(\xi)$ be a polynomial of degree $p$ defined on $I=[-1,1]$ and $q_{i} = q_h(\xi_{i})$, then the expansion of $q_h(\xi)$ in terms of Lagrange
polynomials is given by
            \begin{equation}
                q_h(\xi) := \sum_{i=0}^{p}q_{i}l^{p}_{i}(\xi)\,. \label{eq::nodal_polynomial_expansion}
            \end{equation}
Because the expansion coefficients in \eqref{eq::nodal_polynomial_expansion} are given by the value of $q_h$ in the nodes $\xi_i$,
we refer to this interpolation as a {\em nodal interpolation} and we will denote the Lagrange polynomials in \eqref{eq::lagrange_interpolants} 
by \emph{nodal polynomials}. Using the nodal polynomials we can define another set of basis polynomials, $e^{p}_{i}(\xi)$, as
            \begin{equation}
                e^{p}_{i}(\xi) := - \sum_{k=0}^{i-1}\frac{\mathrm{d}l^{p}_{k}(\xi)}{\mathrm{d}\xi}\,, \qquad i=1,\dots,p\,. \label{eq::histopolant_polynomials_definition}
            \end{equation}
            These polynomials $e^{p}_{i}(\xi)$ have polynomial degree $p-1$ and satisfy,
            \begin{equation}
                 \int_{\xi_{j-1}}^{\xi_{j}}e^{p}_{i}(\xi)\,\mathrm{d}\xi =
                \begin{cases}
                    1 & \mbox{if } i=j \\
                    & \\
                    0 & \mbox{if } i\neq j
                \end{cases}\,, \quad i,j = 1,\dots,p\,. \label{eq::histopolant_polynomials_properties}
            \end{equation}
The proof that the polynomials $e^{p}_{i}(\xi)$ have degree $p-1$ follows directly from the fact that their definition 
\eqref{eq::histopolant_polynomials_definition} involves a linear combination of the derivative of polynomials of degree 
$p$. The proof of \eqref{eq::histopolant_polynomials_properties} results from the properties of $l_{k}^{p}(\xi)$. Using  
\eqref{eq::histopolant_polynomials_definition} the integral of $e^{p}_{i}(\xi)$ becomes
            \[
\int_{\xi_{j-1}}^{\xi_{j}}e^{p}_{i}(\xi)\,\mathrm{d}\xi =  
- \int_{\xi_{j-1}}^{\xi_{j}}\sum_{k=0}^{i-1}\frac{\mathrm{d}l^{p}_{k}(\xi)}{\mathrm{d}\xi} = 
- \sum_{k=0}^{i-1}\int_{\xi_{j-1}}^{\xi_{j}}\frac{\mathrm{d}l^{p}_{k}(\xi)}{\mathrm{d}\xi} = 
- \sum_{k=0}^{i-1} \left(l^{p}_{k}(\xi_{j}) - l^{p}_{k}(\xi_{j-1})\right) = - \sum_{k=0}^{i-1} \left(\delta_{k,j} -\delta_{k,j-1}\right) = 
\delta_{i,j}\,,
            \]
            where $\delta_{i,j}$ is the Kronecker delta.
            For more details see \cite{robidoux-polynomial,Gerritsma11}.

\begin{figure}
    \begin{center}
         \includegraphics[width=0.4\textwidth]{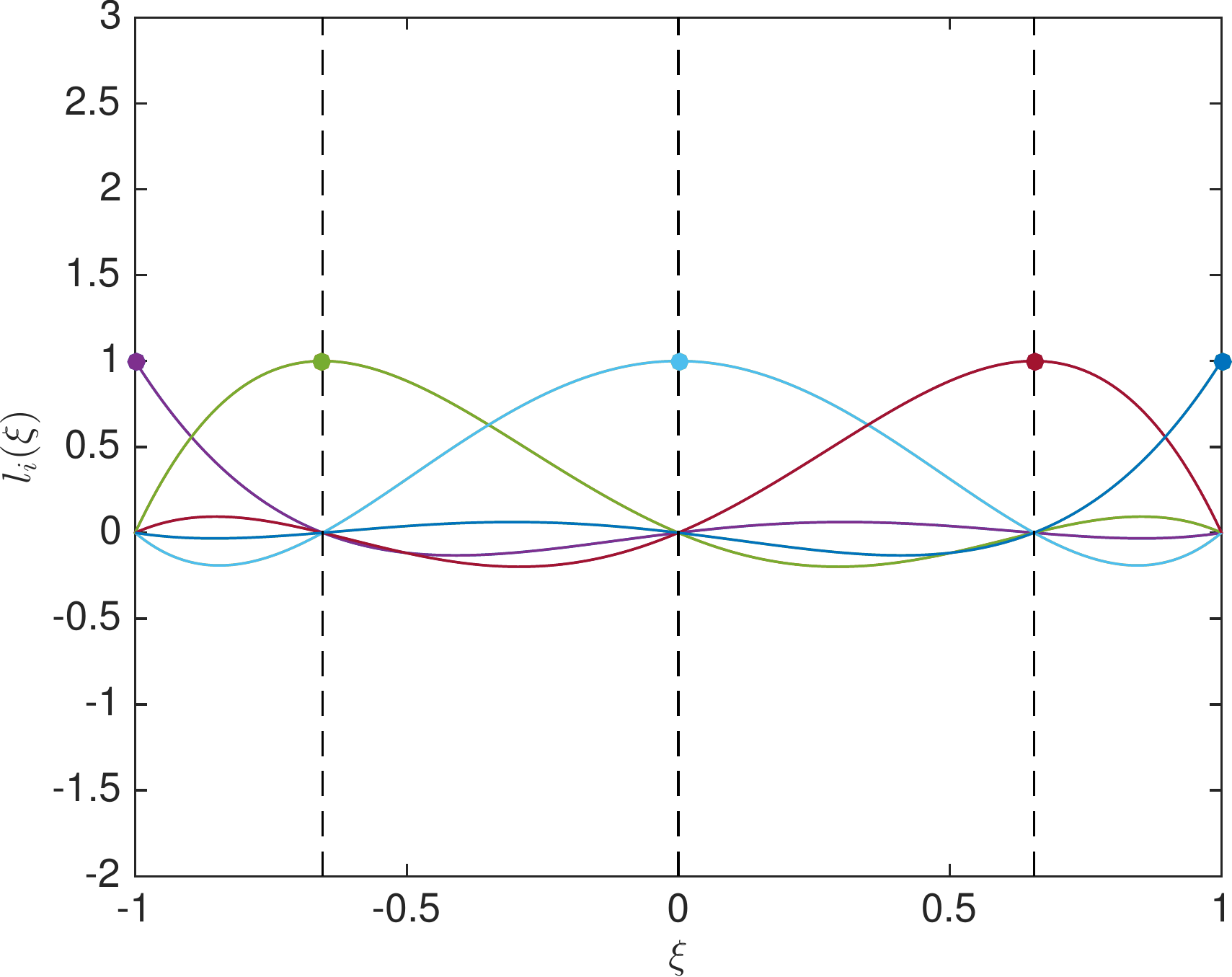} \hspace{1cm}
         \includegraphics[width=0.4\textwidth]{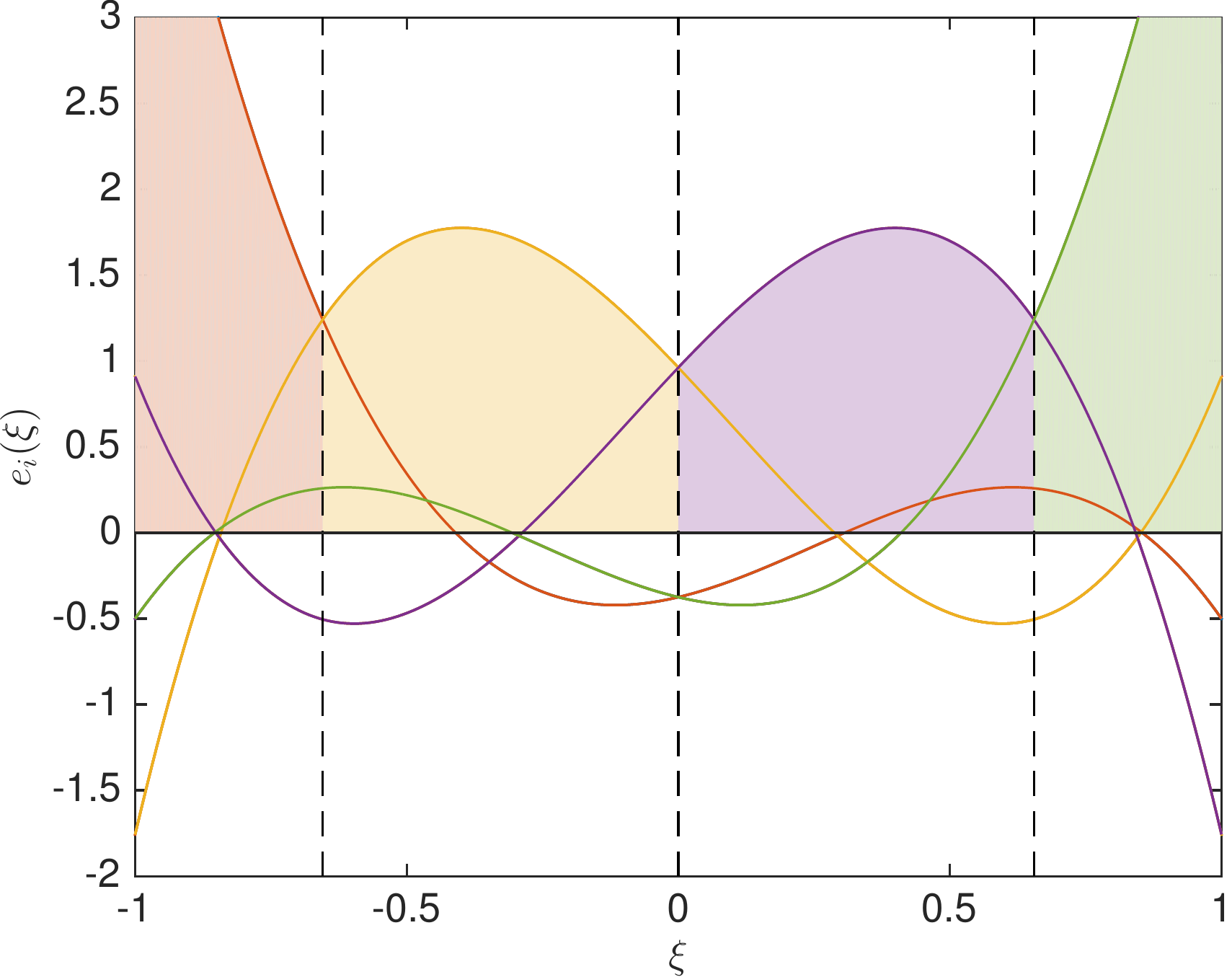}
\caption{Basis polynomials associated to $p=4$. Left: nodal polynomials, the value of the basis polynomial at the 
corresponding node is one and on the other nodes is zero. Right: histopolant polynomials, the integral of the basis 
polynomials over the corresponding shaded area evaluates to one and to zero on the others.}
        \label{fig::basis_polynomials}
    \end{center}
\end{figure}

Let $g_h(\xi)$ be a polynomial of degree $(p-1)$ defined on $I=[-1,1]$ and $g_{i} = \int_{\xi_{i-1}}^{\xi_{i}}g_h(\xi)\,\mathrm{d}\xi$, 
then its expansion in terms of the polynomials $e_i^p(\xi)$ is given by
            \begin{equation}
                g_{h}(\xi) = \sum_{i=1}^{p}g_{i}e^{p}_{i}(\xi)\;.
\label{eq:1D_expansion_in_edge_polynomials}
            \end{equation}
Because the expansion coefficients in \eqref{eq:1D_expansion_in_edge_polynomials} are the integral values of $g_h(\xi)$, 
we denote the polynomials in \eqref{eq::histopolant_polynomials_definition} by \emph{histopolant polynomials} and refer 
to \eqref{eq:1D_expansion_in_edge_polynomials} as {\em histopolation}. It can be shown, \cite{robidoux-polynomial,Gerritsma11}, 
that if $q_h(\xi)$ is expanded in terms of nodal polynomials, as in \eqref{eq::nodal_polynomial_expansion}, then the 
expansion of its derivative $\frac{\mathrm{d}q_h(\xi)}{\mathrm{d}\xi}$ in terms of histopolant polynomials is
\begin{align}
     \left(\frac{\mathrm{d}q_h(\xi)}{\mathrm{d}\xi}\right)_{h} & =
    \sum_{i=1}^{p} \left(\int_{\xi_{i-1}}^{\xi_{i}}\frac{\mathrm{d}q_h(\xi)}{\mathrm{d}\xi}\mathrm{d}\xi\right)e^{p}_{i}(\xi) =
    \sum_{i=1}^{p} \left(q_h(\xi_{i}) - q_h(\xi_{i-1})\right)e^{p}_{i}(\xi) \nonumber \\
 &= \sum_{i=1}^{p} \left(q_{i} - q_{i-1}\right)e^{p}_{i}(\xi) = \sum_{i=1,j=0}^{p}\mathsf{E}^{1,0}_{i,j}q_{j}e^{p}_{i}(\xi)\;,
\end{align}
where $\mathsf{E}^{1,0}_{i,j}$ are the coefficients of the $p\times(p+1)$ matrix $\boldsymbol{\mathsf{E}}^{1,0}$, hereafter
referred to as an \emph{incidence} matrix.
The following identity holds (Commuting property)
\begin{equation}
     \left(\frac{\mathrm{d}q(\xi)}{\mathrm{d}\xi}\right)_{h} = \frac{\mathrm{d}q_{h}(\xi)}{\mathrm{d}\xi}\,.
\end{equation}
For an example of the one-dimensional basis polynomials corresponding to $p=4$, see Fig. \ref{fig::basis_polynomials}.

\subsection{Finite dimensional mimetic function spaces} \label{sec::the_finite_dimensional_basis_functions}

Before discussing the construction of the mixed mimetic spectral elements in two dimensions, we first introduce
the spatial conforming function spaces in which the solution variables will be defined:
\begin{equation}
	W_{h}\subset H(\mathrm{rot},\Omega),\quad U_{h}\subset H(\mathrm{div},\Omega), \quad \mathrm{and} \quad Q_{h}\subset L^{2}(\Omega)\,.
\end{equation}
		
The choice of finite dimensional function spaces determines the properties of the discretization 
\cite{arnold2006finite,arnold2010finite,PG17}, LPG18.
Therefore, the finite dimensional function spaces used in this work are such that when combined form a Hilbert subcomplex
\begin{equation}
	\mathbb{R} \longrightarrow W_{h} \stackrel{\nabla^{\perp}}{\longrightarrow} U_{h} \stackrel{\nabla\cdot}{\longrightarrow} Q_{h}
	\longrightarrow 0\,.\label{eq:hilbert_subcomplex}
\end{equation}
The meaning of this Hilbert subcomplex is that
\begin{equation}
	 \{\nabla^{\perp} \psi_{h}\,|\,\psi_{h}\in W_{h}\} \subset U_{h} \quad \mathrm{and}\quad \{\nabla\cdot\vec{u}_{h}\,|\,\vec{u}_{h}\in U_{h}\} \subseteq Q_{h}\,,
\end{equation}
where $\psi_h$ and $\vec u_h$ are the discrete representations of the stream function and the velocity respectively.
In other words, the \emph{rot} operator must map $W_{h}$ into $U_{h}$ and the \emph{div} operator must map $U_{h}$ onto $Q_{h}$.
		
This discrete subcomplex mimics the 2D Hilbert complex associated to the continuous functional spaces
\begin{equation}
	\mathbb{R} \longrightarrow H(\mathrm{rot},\Omega) \stackrel{\nabla^{\perp}}{\longrightarrow} H(\mathrm{div},\Omega)
	\stackrel{\nabla\cdot}{\longrightarrow} L^{2}(\Omega) \longrightarrow 0\,.\label{eq:hilbert_complex}
\end{equation}
The Hilbert complex is an important structure that is intimately connected to the de Rham complex of differential forms.
Therefore, the construction of a discrete subcomplex is an important requirement to obtain a stable and accurate finite
element discretization
\cite{arnold2010finite,PRHKG14,{bossavit_japan_computational_1,bossavit_japan_computational_2,bossavit_japan_computational_3,bossavit_japan_computational_4,bossavit_japan_computational_5}}.
		
Each of these finite dimensional function spaces $W_{h}$, $U_{h}$, and $Q_{h}$ has an associated finite set of basis functions
$\epsilon_{i}^{\,W}, \vec{\epsilon}_{i}^{\,U}, \epsilon_{i}^{Q}$, such that
\begin{equation}
	W_{h} = \mathrm{span}\{\epsilon_{1}^{\,W},\dots,\epsilon_{d_{W}}^{\,W}\}, \quad U_{h} = \mathrm{span}\{\vec{\epsilon}_{1}^{\,U},\dots,\vec{\epsilon}_{d_{U}}^{\,U}\}, \quad\mathrm{and}\quad Q_{h} = \mathrm{span}\{\epsilon_{1}^{Q},\dots,\epsilon_{d_{Q}}^{Q}\}\,,
\end{equation}
where $d_{W}$, $d_{U}$, and $d_{Q}$ denote the dimension of the discrete function spaces and therefore correspond to the number of
degrees of freedom associated to each of the unknowns.

\subsection{Two dimensional basis functions} \label{sec:two_dimensional_basis_functions}

\textbf{Basis functions for $W_{h}$}
Combining nodal polynomials we can construct the polynomial basis functions for $W_{h}$ on a reference quadrilateral. 
Consider the canonical interval $I=[-1,1]$, the canonical square $\Omega=I\times I\subset \mathbb{R}^{2}$, the nodal 
polynomials \eqref{eq::lagrange_interpolants}, $l^{p}_{i}(\xi)$ of degree $p$, and take $\xi,\eta\in I$. Then a set 
of two-dimensional basis polynomials, $\epsilon^{\,W}_{k}(\xi,\eta; p)$, on $\Omega$ can be constructed as the tensor 
product of the one-dimensional ones
\begin{equation}
    \epsilon^{\,W}_{k}(\xi,\eta; p) := l^{p}_{i}(\xi)\,l^{p}_{j}(\eta)\,, \qquad i,j=0,\dots,p, \quad k = j+1 + i(p+1)\,. \label{eq::2d_nodal_polynomials_definition}
\end{equation}
These polynomials, $\epsilon^{\,W}_{k}(\xi,\eta;p)$, have degree $p$ in each direction and from \eqref{eq::nodal_basis_polynomials} 
it follows that they 
have a value of one at the point $(\xi_i,\eta_j)$ if $k = i(p+1) + j + 1$ and zero otherwise \cite{Gerritsma11,PRHKG14}.
Where, as before, $\xi_{i}$ and $\eta_{i}$ with $i=0,\dots,p$ are the Gauss-Lobatto-Legendre (GLL) nodes. Let $\psi_h(\xi,\eta)$ 
be a polynomial function of degree $p$ in $\xi$ and $\eta$, defined on $\Omega$ and
\begin{equation}
    \psi_{k}^{p} = \psi_h(\xi_{i},\eta_{j}),\qquad \text{with } k = i(p+1) + j+1, \label{eq:2d_nodal_definition_expansion_coefficients}
\end{equation}
then its expansion in terms of these polynomials, $\psi_{h}(\xi,\eta;p)$, is given by
\begin{equation}
    \psi_{h}(\xi,\eta;p) = \sum_{k=1}^{(p+1)^2}\psi_{k}^{p}\,\epsilon_{k}^{\,W}(\xi,\eta;p)\,. \label{eq::2d_nodal_polynomials_expansion}
\end{equation}
For this relation between the expansion coefficients and nodal interpolation we denote the polynomials in \eqref{eq::2d_nodal_polynomials_definition} by
\emph{nodal polynomials}. Therefore we set $W_{h}^{p}:= \mathrm{span}\{\epsilon_{1}^{\,W}(\xi,\eta;p),\dots,\epsilon_{(p+1)^{2}}^{\,W}(\xi,\eta;p)\}$. To
simplify the notation, the explicit reference to the polynomial degree $p$ will be dropped from the function space, the basis functions, and the coefficient
expansion, therefore from here on we will simply use $W_{h}$, $\epsilon_{k}^{\,W}(\xi,\eta)$, and $\psi_{k}$.

\textbf{Basis functions for $U_{h}$}
In a similar fashion we may combine nodal polynomials with histopolant polynomials to construct the polynomial basis 
functions for $U_{h}$ on quadrilaterals. Consider the nodal polynomials \eqref{eq::lagrange_interpolants}, $l^{p}_{i}(\xi)$ 
of degree $p$, the histopolant polynomials \eqref{eq::histopolant_polynomials_definition}, $e^{p}_{i}(\xi)$ of degree $p-1$, 
the canonical square $\Omega=I\times I\subset \mathbb{R}^{2}$, and take $\xi,\eta\in I=[-1,1]$. A set of two-dimensional 
basis polynomials, $\vec{\epsilon}^{\,U}_{k}(\xi,\eta;p)$, can be constructed as the tensor product of the one-dimensional 
basis functions
\begin{equation}
    \vec{\epsilon}^{\,U}_{k}(\xi,\eta;p):=
    \begin{cases}
         l_{i}^{p}(\xi)e_{j}^{p}(\eta)\,\vec{e}_{\xi} & \mbox{if } k \mbox{ odd}, \quad \mbox{with } i=0,\dots,p, \quad j=1,\dots,p, \quad k = 2(ip + j) - 1\;, \\
        & \\
         e_{i}^{p}(\xi)l_{j}^{p}(\eta)\,\vec{e}_{\eta} & \mbox{if } k \mbox{ even}, \quad \mbox{with } i=1,\dots,p, \quad j=0,\dots,p, \quad k = 2((i-1)(p+1) + j + 1)\;.
    \end{cases} \label{eq::edge_polynomials_definition}
\end{equation}
Note that this ordering of degrees of freedom, which alternate between vector components normal to $\vec e_{\xi}$ and those normal to 
$\vec e_{\eta}$ is arbitrary, however it greatly simplifies the implementation, and so will be used in Section \ref{sec::sw}.
These polynomials, $\vec{\epsilon}^{\,U}_{k}(\xi,\eta;p)$, have degree $p$ in $\xi$ and $p-1$ in $\eta$ if $k$ is odd. If $k$ is even, 
then the degree in $\xi$ is $p-1$ and the degree in $\eta$ is $p$.
Let $\vec{u}_h(\xi,\eta;p)$ be a vector valued polynomial function defined on $\Omega$,
then its expansion in terms of these polynomials, $\vec{u}_{h}(\xi,\eta;p)$, is given by
\begin{equation}
    \vec{u}_{h}(\xi,\eta;p) = \sum_{k=1}^{2p(p+1)} u_{k}^{p}\vec{\epsilon}^{\,Q}_{k}(\xi,\eta;p)\,. \label{eq::expansion_edge_polynomials}
\end{equation}
The expansion $\vec{u}_{h}(\xi,\eta;p)$ is a two-dimensional  polynomial edge histopolant (interpolates integral values along lines).
Since the coefficients of this expansion are edge (or flux) integrals, we denote the polynomials in \eqref{eq::edge_polynomials_definition} 
by \emph{edge polynomials}. We set $U_{h} := \mathrm{span}\{\epsilon_{1}^{\,U}(\xi,\eta;p),\dots,\epsilon_{2p(p+1)}^{\,U}(\xi,\eta;p)\}$. 
To simplify the notation, the explicit reference to the polynomial degree $p$ will be dropped from the function space, the basis functions, 
and the expansion coefficients, therefore from here on we will simply use $U_{h}$, $\vec{\epsilon}_{k}^{\,U}(\xi,\eta)$, and $u_{k}$.

\textbf{Basis functions for $Q_{h}$}
Combining histopolant polynomials we can construct the polynomial basis functions for $Q_{h}$ on a quadrilateral. 
Consider the canonical interval $I=[-1,1]$, the canonical square $\Omega=I\times I\subset \mathbb{R}^{2}$, the 
histopolant polynomials \eqref{eq::histopolant_polynomials_definition}, $e^{p}_{i}(\xi)$ of degree $p-1$, and take 
$\xi,\eta\in I$. Then a set of two-dimensional basis polynomials, $\epsilon^{Q}_{k}(\xi,\eta; p)$, can be constructed 
as the tensor product of the one-dimensional ones
\begin{equation}
    \epsilon^{Q}_{k}(\xi,\eta; p) := e^{p}_{i}(\xi)\,e^{p}_{j}(\eta), \qquad i,j=1,\dots,p, \quad k = j + (i-1)p\,. \label{eq::volume_polynomials_definition}
\end{equation}
Where, as before, $\xi_{i}$ and $\eta_{i}$ with $i=0,\dots,p$ are the Gauss-Lobatto-Legendre (GLL) nodes. Let 
$h_h(\xi,\eta)$ be a polynomial function defined on $\Omega$ and 
$h_{k}^{p} = \int_{\xi_{i-1}}^{\xi_{i}}\int_{\eta_{j-1}}^{\eta_{j}}h_h(\xi,\eta)\,\mathrm{d}\xi\mathrm{d}\eta$ 
with $k = j + (i-1)p$, then its expansion in terms of these polynomials, $h_{h}(\xi,\eta;p)$, is given by
\begin{equation}
    h_{h}(\xi,\eta;p) = \sum_{k=1}^{p^{2}}h_{k}^{p}\,\epsilon_{k}^{Q}(\xi,\eta;p)\,. \label{eq::volume_polynomials_expansion}
\end{equation}
For this relation between the expansion coefficients and surface integration we denote the polynomials in 
\eqref{eq::volume_polynomials_definition} by \emph{surface polynomials}. Moreover, these basis polynomials 
satisfy $\epsilon^{Q}_{k}(\xi,\eta;p)\in L^{2}(\Omega)$. Therefore we set 
$Q_{h}^{p}:= \mathrm{span}\{\epsilon_{1}^{Q}(\xi,\eta;p),\dots,\epsilon_{p^{2}}^{Q}(\xi,\eta;p)\}$. To 
simplify the notation, the explicit reference to the polynomial degree $p$ will be dropped from the function 
space, the basis functions, and the coefficient expansion, therefore from here on we will simply use $Q_{h}$, 
$\epsilon_{k}^{Q}(\xi,\eta)$, and $h_{k}$.

\subsection{Properties of the basis functions} \label{sec:properties_of_the_basis_functions}
The first property that can be shown, \cite{Gerritsma11,PRHKG14}, is that if $\psi_h(\xi,\eta) \in W_h$, then
$\nabla^{\perp} \psi_h(\xi,\eta) \in U_h$, where $\psi_h(\xi,\eta)$ is expanded as \eqref{eq::2d_nodal_polynomials_expansion}.

\begin{align}
 \nabla^{\perp}\psi_h(\xi,\eta) 
 & \stackrel{\phantom{\eqref{eq:2d_nodal_definition_expansion_coefficients}}}{=}  \begin{aligned}[t]
          & \sum_{i=0,j=1}^{p}\left(\psi_h(\xi_{i},\eta_{j-1})-\psi_h(\xi_{i},\eta_{j})\right)\vec{\epsilon}_{ip+j}^{\,U}(\xi,\eta) %\nonumber \\
+  \sum_{i=1,j=0}^{p}\left(\psi_h(\xi_{i},\eta_{j}) - \psi_h(\xi_{i-1},\eta_{j})\right)
             \vec{\epsilon}_{(p+i-1)(p+1)+j+1}^{\,U}(\xi,\eta)\nonumber\\
          \end{aligned}\nonumber \\
 &\stackrel{\phantom{\eqref{eq:2d_nodal_definition_expansion_coefficients}}}{=}
\sum_{k=1}^{2p(p+1)} \sum_{j=1}^{(p+1)^2}
\mathsf{E}^{1,0}_{k,j} \psi_{j} \vec{\epsilon}_{k}^{\,U}(\xi,\eta)\,, \label{eq::2d_nodal_basis_expansion_curl}
\end{align}
where $\mathsf{E}_{k,j}^{1,0}$ are the coefficients of the $2p(p+1)\times (p+1)^{2}$ two-dimensional incidence matrix 
$\boldsymbol{\mathsf{E}}^{1,0}$. From \eqref{eq::2d_nodal_basis_expansion_curl} it follows that
\[ \psi_{h}(\xi,\eta) \in W_h \quad \Longrightarrow \quad \nabla^{\perp} \psi_h(\xi,\eta) \in U_h \;,\]
which is the finite-dimensional analogue of
\[ \psi(\xi,\eta) \in H(\mathrm{rot};\Omega) \quad \Longrightarrow \quad \nabla^{\perp} \psi(\xi,\eta) \in H(\mathrm{div};\Omega) \;.\]
Or in fully discrete form, if $\psi_j$ are the expansion coefficients of $\psi_h \in W_h$ with respect to the basis $e_i^W(\xi,\eta)$, then
$\mathsf{E}^{1,0}_{k,j} \psi_{j}$ are the expansion coefficients of $\nabla^{\perp} \psi_h$ in $U_h$ with respect to the basis $\vec{u}_j^U(\xi,\eta)$.
As a special case we have that
\begin{equation}
    \nabla^{\perp} \epsilon^{\,W}_{j} = \sum_{k=1}^{2p(p+1)}\mathsf{E}_{k,j}^{1,0}\vec{\epsilon}^{\,U}_{k}\,, \label{eq:hilbert_subcomplex_basis_W}
\end{equation}
and therefore $\nabla^{\perp}\vec{\epsilon}^{\,W}_{j}\in U_{h}$, with $j=1,\dots,(p+1)^{2}$, and these basis functions satisfy \eqref{eq:hilbert_subcomplex}.

The second property that can be shown, \cite{Gerritsma11,PRHKG14}, is that if $\vec{u}_h(\xi,\eta)\in U_h$ 
is expanded in terms of edge polynomials, as in \eqref{eq::expansion_edge_polynomials}, then the expansion of 
$\nabla\cdot\vec{u}_h$ in terms of the surface polynomials, \eqref{eq::volume_polynomials_definition}, is
\begin{align}
    \nabla\cdot\vec{u}_h(\xi,\eta) &
      =\begin{aligned}[t]
             &\sum_{i,j=1}^{p}\left(\int_{\eta_{j-1}}^{\eta_{i}}\vec{u}_h(\xi_{i},\eta)\cdot\vec{e}_{\xi}\,\mathrm{d}\eta + \int_{\xi_{i-1}}^{\xi_{i}}\vec{u}_h(\xi,\eta_{j})\cdot\vec{e}_{\eta}\,\mathrm{d}\xi -\int_{\eta_{j-1}}^{\eta_{i}}\vec{u}_h(\xi_{i-1},\eta)\cdot\vec{e}_{\xi}\,\mathrm{d}\eta\right. \nonumber \\
           &\qquad\qquad\left.- \int_{\xi_{i-1}}^{\xi_{i}}\vec{u}_h(\xi,\eta_{j-1})\cdot\vec{e}_{\eta}\,\mathrm{d}\xi\right)\epsilon_{j+(i-1)p}^{Q}(\xi,\eta)
        \end{aligned} \nonumber\\
     &                                                                    =  \sum_{k=1}^{p^{2}} \sum_{j=1}^{2p(p+1)}\mathsf{E}^{2,1}_{k,j} u_{j} \epsilon_{k}^{Q}(\xi,\eta)\,, \label{eq::edge_basis_expansion_curl}
\end{align}
where $\mathsf{E}^{2,1}_{k,j}$ are the coefficients of the $p^{2}\times 2p(p+1)$ two-dimensional incidence matrix $\boldsymbol{\mathsf{E}}^{2,1}$. 
Equation \eqref{eq::edge_basis_expansion_curl} confirms that we have a finite dimensional Hilbert sequence as in \eqref{eq:hilbert_subcomplex}, because
\[ \vec{u}_h(\xi,\eta)  \in U_h \quad \Longrightarrow \quad \nabla \cdot \vec{u}_h(\xi,\eta) \in Q_h \;,\]
which is the finite dimensional analogue of
\[ \vec{u} \in H(\mathrm{div};\Omega) \quad \Longrightarrow \quad \nabla \cdot \vec{u} \in L^2(\Omega)\;.\]
In terms of the expansion coefficients we have: If $u_j$ are the expansion coefficients of $\vec{u}_h \in U_h$ 
with respect to the basis $\vec{\epsilon}_j^U$, then the expansion coefficient of $\nabla \cdot \vec{u}_h \in Q_h$ 
with respect to the basis $\epsilon_k^Q$ are given by $\sum_{j=1}^{2p(p+1)}\mathsf{E}^{2,1}_{k,j} u_{j}$.
As a special case we have that
\begin{equation}
\nabla\cdot\vec{\epsilon}^{\,U}_{j} = \sum_{k=1}^{p^{2}}\mathsf{E}_{k,j}^{2,1}\epsilon^{Q}_{k}\,. \label{eq:hilbert_subcomplex_basis_U}
\end{equation}

As seen before, if $\psi_j$ are the expansion coefficients of $\psi_h \in W_h$, then $\mathsf{E}_{k,j}^{1,0} \psi_j$ are the 
expansion coefficients of $\nabla^{\perp} \psi_h \in U_h$. Then $\mathsf{E}_{i,k}^{2,1}\mathsf{E}_{k,j}^{1,0} \psi_j$ are the 
expansion coefficients of $\nabla \cdot \nabla^{\perp} \psi_h \in Q_h$. Since $\nabla \cdot \nabla^{\perp} \psi_h =0$ for all $\psi_h$ 
and because $\epsilon_k^Q$ forms a basis for $Q_h$, we have

\begin{equation}
\sum_{k=1}^{2p(p+1)}\mathsf{E}_{i,k}^{2,1} \circ \mathsf{E}_{k,j}^{1,0} \equiv 0.
\end{equation}
This is the fully discrete representation of the vector identity $\nabla \cdot \nabla^{\perp} \equiv 0$.

In addition to these point wise \emph{strong form} properties, the method also supports corresponding \emph{weak form}
properties via the application of Galerkin projections. Let $\omega_h\in W_h$ and $\vec u_h\in U_h$ 
be discrete representations of the vorticity and velocity respectively, such that $\omega_h = \nabla^*\times\vec u_h$,
where $\nabla^*\times$ is an approximate weak form of the \emph{curl} operator (as opposed to the exact strong form
of the \emph{rot} operator $\nabla^{\perp}$. Then assuming periodic boundary conditions
the adjoint relation between \emph{rot} and \emph{curl} is given as an inner product over the domain $\Omega$ as

\begin{equation}\label{adjt_rot_curl}
\ip{\epsilon^{\,W}_h}{\omega_h} = \ip{\epsilon^{\,W}_h}{\nabla\times\vec u_h} = -\ip{\nabla^{\perp}\epsilon^{\,W}_h}{\vec u_h}.
\end{equation}
The corresponding fully discrete form of this relation is given as 
\begin{equation}\label{curl_wf}
\sum_{i=1}^{d_W}\ip{\epsilon^{\,W}_j}{\epsilon^{\,W}_i}\omega_i = 
-\sum_{k,l=1}^{d_U}(\mathsf{E}^{1,0}_{k,j})^{\top}\ip{\vec\epsilon^{\,U}_k}{\vec\epsilon^{\,U}_l}u_l.
\end{equation}
Similarly, we also have an adjoint relation between \emph{div} and \emph{grad}. Assuming that $p_h\in Q_h$,
for which $\vec u_h = \nabla^*p_h$, where again $\nabla^*$ is an approximate weak form representation of 
\emph{grad} (as opposed to the strong form representation of \emph{div} $\nabla\cdot$), then

\begin{equation}\label{adjt_div_grad}
\ip{\vec\epsilon^{\,U}_h}{\vec u_h} = \ip{\vec\epsilon^{\,U}_h}{\nabla p_h} = -\ip{\nabla\cdot\vec\epsilon^{\,U}_h}{p_h},
\end{equation}
for which the fully discrete form is given as
\begin{equation}\label{grad_wf}
\sum_{i=1}^{d_U}\ip{\vec\epsilon^{\,U}_j}{\vec\epsilon^{\,U}_i}u_i = 
-\sum_{k,l=1}^{d_Q}(\mathsf{E}^{2,1}_{k,j})^{\top}\ip{\epsilon^{\,Q}_k}{\epsilon^{\,Q}_l}p_l.
\end{equation}
From \eqref{curl_wf} and \eqref{grad_wf} it follows that 
$\sum_{i=1}^{d_W}\ip{\epsilon^{\,W}_j}{\epsilon^{\,W}_i}\omega_i = 
\sum_{k=1}^{d_U}\sum_{l,m=1}^{d_Q}(\mathsf{E}^{1,0}_{k,j})^{\top}(\mathsf{E}^{2,1}_{l,k})^{\top}\ip{\epsilon^{\,Q}_l}{\epsilon^{\,Q}_m}p_m$,
such that there is no projection of $p_h$ onto $\omega_h$ since
\begin{equation}\label{curl_grad_wf}
\sum_{j=1}^{2p(p+1)}(\mathsf{E}^{1,0}_{j,i})^{\top}\circ(\mathsf{E}^{2,1}_{k,j})^{\top}\equiv 0,
\end{equation}
which is the \emph{weak form} discrete representation of the vector identity $\nabla\times\nabla\equiv 0$.

\section{Extension to non-affine geometries}\label{sec::met}

		As mentioned before, this work extends the numerical method presented in LPG18 to non-affine geometries. Specifically, 
the main goal being the solution of the shallow water equations on a sphere. For this reason we must introduce how curved 
geometries are treated and included in the mixed mimetic spectral element method. As seen before, the core of the mixed mimetic 
spectral element method is the introduction of geometric degrees of freedom associated with points, lines, surfaces, etc. In turn, 
the mimetic basis functions are interpolatory polynomials that preserve these degrees of freedom. These two aspects are the key 
ingredients required to produce a discretization in which the differential operators can be exactly represented in a purely 
topological manner by incidence matrices (e.g.: $\boldsymbol{\mathsf{E}}^{1,0}, \boldsymbol{\mathsf{E}}^{2,1}$). Therefore, it 
is natural to expect that any extension to curved geometries must be such that these two ingredients are preserved. The 
contravariant and covariant Piola mappings have been thoroughly discussed in \cite{brezzi1991mixed, monk2003} and more recently 
in \cite{RKL09}. It is well known that these two transformations preserve either the tangential or normal traces of vector fields, 
\cite{RKL09}, and have been extensively used for transforming vector fields in $H(\mathrm{div})$ and $H(\mathrm{curl})$ respectively. 
We will construct transformation rules $\Phi_{0}^{*}$, $\Phi_{1,\mathrm{curl}}^{*}$, $\Phi_{1,\mathrm{div}}^{*}$, and $\Phi_{2}^{*}$ 
for $H(\mathrm{rot}, \Omega)$, $H(\mathrm{curl}, \Omega)$, $H(\mathrm{div}, \Omega)$, and $L^{2}(\Omega)$, respectively. To do so, we 
will require that the geometric degrees of freedom are preserved by the transformation. The subscript $i$ in the transformations 
$\Phi_{i}^{*}$ explicitly establish this relation to the geometric degrees of freedom: (0) nodal degrees of freedom, (1) edge 
degrees of freedom, and (2) surface degrees of freedom. We show that for vector fields this leads to the covariant and contravariant 
Piola transformations.
		 
		Additionally, we will also introduce the commuting relation between the transformations and the differential operators and 
the invariance of the incidence matrices (discrete differential operators). To finalize this discussion on curved geometries, we 
will explicitly show how the curved geometry affects the inner products used throughout this work and how they may be generalized 
with respect to the work presented in LPG18.
		
	\subsection{Transformation rules}
		\subsubsection{Scalar fields}\label{sec:transformation_rules_scalar_fields}
			Consider the two-dimensional manifolds $\mathcal{M}$ and $\tilde{\mathcal{M}}$, and the nondegenerate mapping $\Phi: \mathcal{M}\mapsto \tilde{\mathcal{M}}$ with Jacobian $\boldsymbol{\mathsf{J}}$ and Jacobian determinant $J$. For $\tilde{f}\in H(\mathrm{rot},\tilde{\mathcal{M}})$ we introduce the mapping $\Phi^{*}_{0}: H(\mathrm{rot},\tilde{\mathcal{M}})\mapsto H(\mathrm{rot},\mathcal{M})$
			\begin{equation}
				\Phi^{*}_{0}\left[\tilde{f}\right] := \tilde{f}\circ\Phi\,, \label{eq:pullback_0_forms}
			\end{equation}
			and we say that $f := \tilde{f}\circ\Phi$ with $f\in H(\mathrm{rot},\mathcal{M})$. Moreover, the inverse mapping $\left(\Phi^{*}_{0}\right)^{-1}:H(\mathrm{rot},\mathcal{M})\mapsto H(\mathrm{rot},\tilde{\mathcal{M}})$ is
			\begin{equation}
				\left(\Phi^{*}_{0}\right)^{-1}\left[f\right] := f\circ\Phi^{-1} = \tilde{f}\,. \label{eq:pullback_inverse_0_forms}
			\end{equation}
			
			The significance of this transformation is that if $(\xi^{1}, \xi^{2}) \in \mathcal{M}$ and $\Phi(\xi^{1}, \xi^{2}) = (x^{1}, x^{2}) \in \tilde{\mathcal{M}}$, then
			\begin{equation}
				f(\xi^{1},\xi^{2}) = \left(\Phi^{*}_{0}\right)^{-1}\left[f\right] (x^{1}, x^{2}), \qquad\mathrm{and}\qquad \tilde{f}(x^{1},x^{2}) = \Phi^{*}_{0}\left[\tilde{f}\right] (\xi^{1}, \xi^{2})\,,
			\end{equation}
			that is: point evaluations are preserved under this transformation. Therefore, the geometric degrees of freedom used to discretize functions in $H(\mathrm{rot},\mathcal{M})$ are invariant under this transformation.
			
			We now introduce another transformation, specifically constructed for scalar fields $\tilde{f}\in L^{2}(\tilde{\mathcal{M}})$. 
As seen before, the geometric degrees of freedom associated to these scalar fields are volume integrations (since we are considering only 
two dimensional manifolds, volumes become surfaces). It is well known, see for example \cite[Section~12.7]{apostol1969vol2}, that the 
surface integral of a scalar field $\tilde{f}\in L^{2}(\tilde{\mathcal{M}})$ over a manifold $\tilde{\mathcal{M}}$ is related to the 
integral over the manifold $\mathcal{M}$ in the following way
			\begin{equation}
				\int_{\tilde{\mathcal{M}}} \tilde{f} \,\mathrm{d}\tilde{\mathcal{M}} = \int_{\mathcal{M}}\left(\tilde{f}\circ\Phi\right) \, J\,\mathrm{d}\mathcal{M}\,.
			\end{equation}
			We can then introduce the transformation $\Phi^{*}_{2}: L^{2}(\tilde{\mathcal{M}})\mapsto L^{2}(\mathcal{M})$ as
			\begin{equation}
				\Phi^{*}_{2}\left[\tilde{f}\right] := \left(\tilde{f}\circ\Phi\right)\,J\,, \label{eq:pullback_2_forms}
			\end{equation}
			and we say that $f := \left(\tilde{f}\circ\Phi\right)\,J$ with $f\in L^{2}(\mathcal{M})$. In a similar way as before, the inverse mapping $\left(\Phi_{2}^{*}\right)^{-1}:L^{2}(\mathcal{M})\mapsto L^{2}(\tilde{\mathcal{M}})$ is
			\begin{equation}
				\left(\Phi^{*}_{2}\right)^{-1}\left[f\right] := \left(f\circ\Phi^{-1}\right)\,\frac{1}{J} = \tilde{f}\,. \label{eq:pullback_inverse_2_forms}
			\end{equation}
			
			The relevance of this transformation, analogously to the previous one, lies in the fact that for submanifolds $\mathcal{N}\subseteq\mathcal{M}$ and $\Phi(\mathcal{N}) = \tilde{\mathcal{N}}\subseteq\tilde{\mathcal{M}}$ we have
			\begin{equation}
				\int_{\mathcal{N}} f\,\mathrm{d}\mathcal{N} = \int_{\tilde{\mathcal{N}}} \left(\Phi^{*}_{2}\right)^{-1}\left[f\right]\,\mathrm{d}\tilde{\mathcal{N}} \qquad \mathrm{and} \qquad \int_{\tilde{\mathcal{N}}} \tilde{f}\,\mathrm{d}\tilde{\mathcal{N}} = \int_{\mathcal{N}} \Phi^{*}_{2}\left[\tilde{f}\right]\,\mathrm{d}\mathcal{N}\,. \label{eq:pullback_2_forms_integral_invariance}
			\end{equation}
			Under this transformation, surface integrals are preserved. This means that the geometric degrees of freedom associated 
with the discretization of functions in $L^{2}(\mathcal{M})$ are invariant under this transformation.
			
		\subsubsection{Vector fields}\label{sec:transformation_rules_vector_fields}
			As mentioned in the start of this section, the transformation rules for vector fields are the covariant and contravariant 
Piola transformations. We will briefly show that these transformations derive directly from the same ideas used to derive the 
transformations for the scalar fields, i.e. invariance of the geometric degrees of freedom. In the case of vector fields the geometric 
degrees of freedom are line integrals for $H(\mathrm{curl},\mathcal{M})$, and flux integrals for $H(\mathrm{div},\mathcal{M})$ (since 
here we consider only two dimensional manifolds, the flux integrals become normal line integrals).
			
			As usual, see for example \cite[Section~10.2]{apostol1969vol2}, the line integral of a vector field 
$\tilde{\vec{v}}\in H(\mathrm{curl},\tilde{\mathcal{M}})$ along a line segment $\tilde{\gamma}:I\subset\mathbb{R}\mapsto\tilde{\mathcal{M}}$ is
			\begin{equation}
				\int_{\tilde{\gamma}}\left(\tilde{\vec{v}}\circ\tilde{\gamma}\right)\cdot \tilde{\vec{t}}\,\mathrm{d}\tilde{\gamma} := 
\int_{I}\left(\tilde{\vec{v}}\circ\tilde{\gamma}\right)\cdot \frac{\vec{\mathrm{d}\tilde{\gamma}}}{\mathrm{d}s}\,\mathrm{d}s\,,\label{eq:line_integral}
			\end{equation}
where $\tilde{\vec t}$ is a tangent vector and $s$ is a parametric coordinate along the line.
			In a similar manner, the line integral of  a vector field $\vec{v}\in H(\mathrm{curl},\mathcal{M})$ along a line segment $\gamma:I\subset\mathbb{R}\mapsto\mathcal{M}$ is
			\begin{equation}
				\int_{\gamma}\left(\vec{v}\circ\gamma\right)\cdot \vec{t}\,\mathrm{d}\gamma := \int_{I}\left(\vec{v}\circ\gamma\right)\cdot \frac{\vec{\mathrm{d}\gamma}}{\mathrm{d}s}\,\mathrm{d}s\,. 
			\end{equation}
			If we consider the case $\tilde{\gamma} = \Phi[\gamma]$, then we may rewrite \eqref{eq:line_integral} as
			\begin{equation}
				\int_{I}\left(\tilde{\vec{v}}\circ\tilde{\gamma}\right)\cdot \frac{\vec{\mathrm{d}\tilde{\gamma}}}{\mathrm{d}s}\,\mathrm{d}s = \int_{I}\left(\tilde{\vec{v}}\circ\Phi\circ\gamma\right)\cdot \frac{\vec{\mathrm{d}\left(\Phi\circ\gamma\right)}}{\mathrm{d}s}\,\mathrm{d}s = \int_{I}\left(\tilde{\vec{v}}\circ\Phi\circ\gamma\right)\cdot \boldsymbol{\mathsf{J}}\frac{\vec{\mathrm{d}\gamma}}{\mathrm{d}s}\,\mathrm{d}s = \int_{I}\boldsymbol{\mathsf{J}}^{\top}\left(\tilde{\vec{v}}\circ\Phi\circ\gamma\right)\cdot \frac{\vec{\mathrm{d}\gamma}}{\mathrm{d}s}\,\mathrm{d}s\,, \label{eq:line_integral_2}
			\end{equation}
where the Jacobian $\boldsymbol{\mathsf{J}}$ is defined such that $J_{i,j} = \frac{\partial\Phi^i}{\partial\xi^j}$. Now, if we want to construct a 
transformation $\Phi^{*}_{1,\mathrm{curl}}: H(\mathrm{curl},\tilde{\mathcal{M}})\mapsto H(\mathrm{curl},\mathcal{M})$ that 
preserves line integrals, then for $\tilde{\gamma} = \Phi[\gamma]$ and $\vec{v} = \Phi^{*}_{1,\mathrm{curl}}\left[\tilde{\vec{v}}\right]$ we 
must to satisfy
			\begin{equation}
				\int_{I}\left(\tilde{\vec{v}}\circ\tilde{\gamma}\right)\cdot \frac{\vec{\mathrm{d}\tilde{\gamma}}}{\mathrm{d}s}\,\mathrm{d}s = \int_{I}\left(\vec{v}\circ\gamma\right)\cdot \frac{\vec{\mathrm{d}\gamma}}{\mathrm{d}s}\,\mathrm{d}s\,. \label{eq:line_integral_3}
			\end{equation}			
			Therefore, combining \eqref{eq:line_integral_2}  with \eqref{eq:line_integral_3} we obtain
			\begin{equation}
				\vec{v} = \boldsymbol{\mathsf{J}}^{\top}\left(\tilde{\vec{v}}\circ\Phi\right) := \Phi^{*}_{1,\mathrm{curl}}\left[\tilde{\vec{v}}\right]\,. \label{eq:pullback_inner_1_forms}
			\end{equation}
			The inverse mapping  $\left(\Phi^{*}_{1,\mathrm{curl}}\right)^{-1}: H(\mathrm{curl},\mathcal{M})\mapsto H(\mathrm{curl},\tilde{\mathcal{M}})$ follows directly
			\begin{equation}
				\left(\Phi^{*}_{1,\mathrm{curl}}\right)^{-1}\left[\vec{v}\right] := \boldsymbol{\mathsf{J}}^{-\top}\left(\vec{v}\circ\Phi^{-1}\right) = \tilde{\vec{v}}\,. \label{eq:pullback_inverse_inner_1_forms}
			\end{equation}
			The inverse transformation $\left(\Phi^{*}_{1,\mathrm{curl}}\right)^{-1}$ is the \emph{covariant} Piola transformation. 
As we have seen, by construction, this transformation preserves line integrals. Therefore, the geometric degrees of freedom associated 
to the discretization of vector fields in $H(\mathrm{curl},\mathcal{M})$ are kept constant under this transformation.
			
			The construction of the transformation for $H(\mathrm{div},\tilde{\mathcal{M}})$ that preserves normal line integrals 
(fluxes in two dimensions) is almost identical to the one we just derived for $H(\mathrm{curl},\tilde{\mathcal{M}})$. The only 
difference being the use of the normal to the curve, $\tilde{\vec n}$, instead of the tangent. Therefore \eqref{eq:line_integral} becomes
			\begin{equation}
				\int_{\tilde{\gamma}}\left(\tilde{\vec{v}}\circ\tilde{\gamma}\right)\cdot \tilde{\vec{n}}\,\mathrm{d}\tilde{\gamma} := \int_{I}\left(\tilde{\vec{v}}\circ\tilde{\gamma}\right)\cdot \boldsymbol{\mathsf{H}}\frac{\vec{\mathrm{d}\tilde{\gamma}}}{\mathrm{d}s}\,\mathrm{d}s\,, \qquad \mathrm{with}\qquad \boldsymbol{\mathsf{H}} := \pm\left[\begin{array}{cc}0 & -1 \\ 1 & 0\end{array}\right]\,,
			\end{equation}
			depending on the chosen orientation, either the plus or the minus sign should be used. If we now follow exactly the same procedure as before, but taking into consideration the additional matrix $\boldsymbol{\mathsf{H}}$, we obtain the following transformation rule $\Phi^{*}_{1,\mathrm{div}}: H(\mathrm{div},\tilde{\mathcal{M}})\mapsto H(\mathrm{div},\mathcal{M})$ that preserves fluxes
			\begin{equation}
				\Phi^{*}_{1,\mathrm{div}}\left[\tilde{\vec{v}}\right] := J \boldsymbol{\mathsf{J}}^{-1}\left(\tilde{\vec{v}}\circ\Phi\right) \,, \label{eq:pullback_outer_1_forms}
			\end{equation}
			and we say that $\vec{v}:= J \boldsymbol{\mathsf{J}}^{-1}\left(\tilde{\vec{v}}\circ\Phi\right)$ with $\vec{v}\in H(\mathrm{div}, \mathcal{M})$. The inverse transform $\left(\Phi^{*}_{1,\mathrm{div}}\right)^{-1}: H(\mathrm{div},\mathcal{M})\mapsto H(\mathrm{div},\tilde{\mathcal{M}})$ can be easily computed and is given by
			\begin{equation}
				\left(\Phi^{*}_{1,\mathrm{div}}\right)^{-1}\left[\vec{v}\right] := \frac{1}{J}\,\boldsymbol{\mathsf{J}}\left(\vec{v}\circ\Phi^{-1}\right) = \tilde{\vec{v}} \,. \label{eq:pullback_inverse_outer_1_forms}
			\end{equation}
			This inverse transformation $\left(\Phi^{*}_{1,\mathrm{div}}\right)^{-1}$ is the \emph{contravariant} Piola transformation. 
Since flux integrals are preserved by this transformation, also the geometric degrees of freedom associated to the discretization of 
vector fields in $H(\mathrm{div},\mathcal{M})$ remain invariant.

	\subsection{Commuting relations for the differential operators}
			Having introduced the transformation rules for all types of physical field quantities present in this work
\footnote{Although not directly used in this work, we have also presented the transformation rule for the 
$H(\mathrm{curl},\mathcal{M})$ space, for completeness.}, the next natural step is to discuss the relation between the 
differential operators and these transformations. Due to their length, we will not present full demonstrations. Instead, 
the results will be given, followed by a brief outline for the cases used in this work: $\nabla^{\perp}$ and $\nabla\cdot$.
			
			The key result is that the differential operators $\nabla^{\perp}$, $\nabla\times$, $\nabla\cdot$ commute 
with the transformations introduced in the previous section, that is
			\begin{equation}
				\Phi^{*}_{1,\mathrm{div}}\left[\tilde{\nabla}^{\perp}\tilde{f}\right] = \nabla^{\perp}\Phi^{*}_{0}\left[\tilde{f}\right], \qquad \Phi^{*}_{2}\left[\tilde{\nabla}\times\tilde{\vec{v}}\right] = \nabla\times\Phi^{*}_{1,\mathrm{curl}}\left[\tilde{\vec{v}}\right], \qquad\mathrm{and}\qquad \Phi^{*}_{2}\left[\tilde{\nabla}\cdot\tilde{\vec{u}}\right] = \nabla\cdot\Phi^{*}_{1,\mathrm{div}}\left[\tilde{\vec{u}}\right]\,, \label{eq:commuting_relations_nabla}
			\end{equation}
			with $\tilde{f}\in H(\mathrm{rot},\tilde{\mathcal{M}})$, $\tilde{\vec{v}}\in H(\mathrm{curl},\tilde{\mathcal{M}})$, 
and $\tilde{\vec{u}}\in H(\mathrm{div},\tilde{\mathcal{M}})$ . Where as usual 
$\tilde{\nabla}^{\perp}\tilde{f} = -\frac{\partial\tilde{f}}{\partial x^{2}}\vec{e}_{1} + \frac{\partial\tilde{f}}{\partial x^{1}}\vec{e}_{2}$, 
$\tilde{\nabla}\times\tilde{\vec{v}} = \frac{\partial\tilde{v}^{2}}{\partial x^{1}} - \frac{\partial\tilde{v}^{1}}{\partial x^{2}}$, 
and $\tilde{\nabla}\cdot\tilde{\vec{u}} = \frac{\partial\tilde{u}^{1}}{\partial x^{1}} + \frac{\partial\tilde{u}^{2}}{\partial x^{2}}$. 
The differential operators in $\mathcal{M}$ are identical, with the exception that $x^{i}$ is replaced by $\xi^{i}$, with $i=\{1,2\}$. 
An important aspect to highlight here is the pairs of transformations used in the commuting relations. These pairs derive directly 
from the Hilbert subcomplex \eqref{eq:hilbert_subcomplex} (and an analogous one for $\nabla$ and $\nabla\times$). For example, the 
pair  $\Phi^{*}_{0}$ and $\Phi^{*}_{1,\mathrm{div}}$ is used because if $\tilde{f}\in H(\mathrm{rot},\tilde{\mathcal{M}})$ then 
$\tilde{\nabla}^{\perp}\tilde{f}\in H(\mathrm{div},\tilde{\mathcal{M}})$, as is the case in \eqref{eq:commuting_relations_nabla}.
		
		\subsubsection{Commuting relation for $\nabla^{\perp}$}
			To show the commuting relation for $\nabla^{\perp}$ we start with the commuting relation itself and explicitly substitute 
the transformations $\Phi^{*}_{0}$ and $\Phi^{*}_{1,\mathrm{div}}$, respectively  \eqref{eq:pullback_0_forms} and  \eqref{eq:pullback_outer_1_forms},
			\begin{equation}
				J \boldsymbol{\mathsf{J}}^{-1}\left(\left(\tilde{\nabla}^{\perp}\tilde{f}\right)\circ\Phi\right) := \Phi^{*}_{1,\mathrm{div}}\left[\tilde{\nabla}^{\perp}\tilde{f}\right] = \nabla^{\perp}\Phi^{*}_{0}\left[\tilde{f}\right] =: \nabla^{\perp} \left(\tilde{f}\circ\Phi\right)\,.
			\end{equation}
			After substitution of the inverse matrix Jacobian and expanding the derivatives on the right hand side using the chain rule we get the following identity
			\begin{equation}
				\left[
					\begin{array}{*2{>{\displaystyle}c}p{5cm}}
						\frac{\partial\Phi^{2}}{\partial \xi^{2}} & -\frac{\partial\Phi^{1}}{\partial \xi^{2}} \\
						-\frac{\partial\Phi^{2}}{\partial \xi^{1}} & \frac{\partial\Phi^{1}}{\partial \xi^{1}}
					\end{array}
				\right]
				\left[
					\begin{array}{*1{>{\displaystyle}c}p{5cm}}
						-\frac{\partial\tilde{f}}{\partial x^{2}} \\
						\frac{\partial\tilde{f}}{\partial x^{1}}
					\end{array}
				\right]
				=
				\left[
					\begin{array}{*2{>{\displaystyle}c}p{5cm}}
						\frac{\partial\Phi^{2}}{\partial \xi^{2}} & -\frac{\partial\Phi^{1}}{\partial \xi^{2}} \\
						-\frac{\partial\Phi^{2}}{\partial \xi^{1}} & \frac{\partial\Phi^{1}}{\partial \xi^{1}}
					\end{array}
				\right]
				\left[
					\begin{array}{*1{>{\displaystyle}c}p{5cm}}
						-\frac{\partial\tilde{f}}{\partial x^{2}} \\
						\frac{\partial\tilde{f}}{\partial x^{1}}
					\end{array}
				\right]\,,
			\end{equation}
			thus proving the commuting relation.
			
		\subsubsection{Commuting relation for $\nabla\cdot$}
			The proof for the commuting relation of $\nabla\cdot$ follows the same steps. The first step is the commuting relation itself with the explicit expressions for the transformations $\Phi_{2}^{*}$ and $\Phi_{1\mathrm{div}}^{*}$, respectively \eqref{eq:pullback_2_forms} and  \eqref{eq:pullback_outer_1_forms},
			\begin{equation}
				J\left(\tilde{\nabla}\cdot\tilde{\vec{u}}\right)\circ\Phi := \Phi^{*}_{2}\left[\tilde{\nabla}\cdot\tilde{\vec{u}}\right] = \nabla\cdot\Phi^{*}_{1,\mathrm{div}}\left[\tilde{\vec{u}}\right] =: \nabla\cdot \left(J \boldsymbol{\mathsf{J}}^{-1}\left(\tilde{\vec{u}}\circ\Phi\right)\right)\,.
			\end{equation}
			If we expand both the right hand side and the left hand side expression we get
			\begin{equation}
				J\left(\frac{\partial\tilde{u}^{1}}{\partial x^{1}}\circ\Phi + \frac{\partial\tilde{u}^{2}}{\partial x^{2}}\circ\Phi\right) = \frac{\partial}{\partial\xi^{1}}\left(\frac{\partial\Phi^{2}}{\partial \xi^{2}} (\tilde{u}^{1}\circ\Phi) - \frac{\partial\Phi^{1}}{\partial \xi^{2}} (\tilde{u}^{2}\circ\Phi)\right) + \frac{\partial}{\partial\xi^{2}}\left(-\frac{\partial\Phi^{2}}{\partial \xi^{1}} (\tilde{u}^{1}\circ\Phi) + \frac{\partial\Phi^{1}}{\partial \xi^{1}} (\tilde{u}^{2}\circ\Phi)\right)\,.
			\end{equation}
			Expanding the right hand side using the chain rule and after some manipulations we obtain the identity
			\begin{equation}
				J\left(\frac{\partial\tilde{u}^{1}}{\partial x^{1}}\circ\Phi + \frac{\partial\tilde{u}^{2}}{\partial x^{2}}\circ\Phi\right) = J\left(\frac{\partial\tilde{u}^{1}}{\partial x^{1}}\circ\Phi + \frac{\partial\tilde{u}^{2}}{\partial x^{2}}\circ\Phi\right)\,,
			\end{equation}			
			finishing the proof.
		
	\subsection{Inner products in $H(\mathrm{rot},\Omega)$, $H(\mathrm{div},\Omega)$, and $L^{2}(\Omega)$}
		Another important point to discuss in the context of curved geometries is how to compute the inner products that appear in this work, namely
		\begin{equation}
			\langle \tilde{f}, \tilde{g}\rangle_{\Omega} := \int_{\Omega}\tilde{f}\tilde{g}\,\mathrm{d}\Omega, \quad \langle \tilde{\vec{v}}, \tilde{\vec{u}}\rangle_{\Omega} := \int_{\Omega}\tilde{\vec{v}}\cdot\tilde{\vec{u}}\,\mathrm{d}\Omega, \quad\mathrm{and}\quad \langle \tilde{q}, \tilde{a}\rangle_{\Omega} := \int_{\Omega}\tilde{q}\tilde{a}\,\mathrm{d}\Omega\,, \label{eq:continuous_inner_product_curved_geometries}
		\end{equation}
		with $\tilde{f},\tilde{g}\in H(\mathrm{rot},\Omega)$, $\tilde{\vec{v}},\tilde{\vec{u}}\in H(\mathrm{div},\Omega)$, 
and $\tilde{q},\tilde{a}\in L^{2}(\Omega)$. It is straightforward to note that all integrands that appear in 
\eqref{eq:continuous_inner_product_curved_geometries} are scalar fields. Therefore, since the goal is to compute these 
scalar integral quantities we must use a transformation rule that preserves scalar surface integrals, that is 
\eqref{eq:pullback_2_forms}. To apply \eqref{eq:pullback_2_forms} there must exist a nondegenerate mapping 
$\Phi: \Omega_{0}\mapsto\Omega$, with $\Omega_{0} = [-1,1]\times[-1,1]$ since here we consider only the two-dimensional 
case. Under these conditions \eqref{eq:continuous_inner_product_curved_geometries} may be transformed into
		\begin{equation}
			\langle \tilde{f}, \tilde{g}\rangle_{\Omega} := \int_{\Omega_{0}}\left(\tilde{f}\circ\Phi\right)\left(\tilde{g}\circ\Phi\right)\, J\mathrm{d}\Omega_{0}, \quad \langle \tilde{\vec{v}}, \tilde{\vec{u}}\rangle_{\Omega} := \int_{\Omega_{0}}\left(\tilde{\vec{v}}\circ\Phi\right)\cdot\left(\tilde{\vec{u}}\circ\Phi\right)\,J\mathrm{d}\Omega_{0}, \quad\mathrm{and}\quad \langle \tilde{q}, \tilde{a}\rangle_{\Omega} := \int_{\Omega_{0}}\left(\tilde{q}\circ\Phi\right)\left(\tilde{a}\circ\Phi\right)\,J\mathrm{d}\Omega_{0}\,, \label{eq:continuous_inner_product_curved_geometries_2}
		\end{equation}
		which can now be straightforwardly integrated.
		
	\subsection{Inner products in $\tilde{W}_{h}\subset H(\mathrm{rot},\Omega)$, $\tilde{U}_{h}\subset H(\mathrm{div},\Omega)$, and $\tilde{Q}_{h}\subset L^{2}(\Omega)$}
		As is usual in all finite element formulations, the domain $\Omega$ is partitioned into a set of $K$ 
non-overlapping subdomains $\Omega_{i}$ with $i=1,\dots,K$ such that $\Omega = \bigcup_{i=1}^{K}\Omega_{i}$. 
Moreover, we consider that there exist $K$ nondegenerate mappings $\Phi_{i}:\Omega_{0}\mapsto\Omega_{i}$.  
Therefore we may focus on one single subdomain $\Omega_{i}$ since the results apply to all subdomains. To 
simplify the notation, on what follows we will suppress the index $i$ and use $\Omega$ and $\Phi$ instead of 
$\Omega_{i}$ and $\Phi_{i}$.
		
		\subsubsection{$\tilde{W}_{h}\subset H(\mathrm{rot},\Omega)$}		
			In Section~\ref{sec:two_dimensional_basis_functions} we introduced the basis functions for the space $W_{h}\subset H(\mathrm{rot},\Omega_{0})$ such that 
			\[
				W_{h} := \mathrm{span}\{\epsilon_{1}^{W}(\xi,\eta), \dots, \epsilon^{W}_{(p+1)^{2}}(\xi,\eta)\}\,,
			\]
			where $p$ is the polynomial degree of the basis functions, and $(\xi,\eta)\in\Omega_{0}=[-1,1]\times[-1,1]$. We also saw that any function $f_{h}\in W_{h}$ could be written as a linear combination of the basis functions, \eqref{eq::2d_nodal_polynomials_expansion},
			\begin{equation}
				f_{h} = \sum_{i=1}^{(p+1)^{2}}f_{i}\,\epsilon_{i}^{W}\,.
			\end{equation}
			
			If we apply the transformation rule $\left(\Phi^{*}_{0}\right)^{-1}$, \eqref{eq:pullback_inverse_0_forms}, to $f_{h} \in W_{h}$ we obtain $\tilde{f}_{h}\in \tilde{W}_{h}$, its polynomial expansion in the physical domain $\Omega$
			\begin{equation}
				\tilde{f}_{h} := \sum_{i=1}^{(p+1)^{2}}f_{i}\,\left(\Phi^{*}_{0}\right)^{-1}\left[\epsilon_{i}^{W}\right] := \sum_{i=1}^{(p+1)^{2}}f_{i}\,\left(\epsilon_{i}^{W}\circ\Phi^{-1}\right) \, := \sum_{i=1}^{(p+1)^{2}}f_{i}\,\tilde{\epsilon}_{i}^{W}\,, \label{eq:0_form_expansion_curved_domain}
			\end{equation}
			where we have defined the basis functions $\tilde{\epsilon}_{i}^{W}$ of the space $\tilde{W}_{h}$ as 
			\begin{equation}
				\tilde{\epsilon}_{i}^{W} := \left(\epsilon_{i}^{W}\circ\Phi^{-1}\right)\,, \qquad i=1,\dots,(1+p)^{2}\,. \label{eq:basis_functions_0_omega}
			\end{equation}
			Therefore, the inner product between two scalar functions $\tilde{f}_{h}, \tilde{g}_{h}\in \tilde{W}_{h}$ is
			\begin{equation}
				\langle\tilde{f}_{h}, \tilde{g}_{h}\rangle_{\Omega} \stackrel{\eqref{eq:0_form_expansion_curved_domain}}{=} \langle\sum_{i=1}^{(p+1)^{2}}f_{i}\,\left(\epsilon_{i}^{W}\circ\Phi^{-1}\right), \sum_{j=1}^{(p+1)^{2}}g_{j}\,\left(\epsilon_{j}^{W}\circ\Phi^{-1}\right)\rangle_{\Omega} \stackrel{\eqref{eq:continuous_inner_product_curved_geometries}}{=} \sum_{i,j}^{(p+1)^{2}}f_{i}\,g_{j}\int_{\Omega}\left(\epsilon_{i}^{W}\circ\Phi^{-1}\right) \left(\epsilon_{j}^{W}\circ\Phi^{-1}\right)\,\mathrm{d}\Omega\,.
			\end{equation}
			As before, we may now apply the transformation $\Phi^{*}_{2}$, \eqref{eq:pullback_2_forms}, to the integrand in order to transform the integral into an integral in $\Omega_{0}$ without changing the value of the integral
			\begin{equation}
				\langle\tilde{f}_{h}, \tilde{g}_{h}\rangle_{\Omega} = \sum_{i,j}^{(p+1)^{2}}f_{i}\,g_{j}\int_{\Omega}\left(\epsilon_{i}^{W}\circ\Phi^{-1}\right) \left(\epsilon_{j}^{W}\circ\Phi^{-1}\right)\,\mathrm{d}\Omega \stackrel{\eqref{eq:pullback_2_forms_integral_invariance}}{=} \sum_{i,j}^{(p+1)^{2}}f_{i}\,g_{j}\int_{\Omega_{0}}\epsilon_{i}^{W} \epsilon_{j}^{W}\,J\mathrm{d}\Omega\,.
			\end{equation}
			To simplify the notation, we define an inner product in $W_{h}$, $\langle\cdot,\cdot\rangle : W_{h}\times W_{h} \mapsto \mathbb{R}$ as
			\begin{equation}\label{eq:inner_product_Wh}
				\langle f_{h}, g_{h}\rangle_{\Omega_{0}} := \int_{\Omega_{0}} f_{h}g_{h} \, J\mathrm{d}\Omega_{0}\,, \quad \text{with}\quad f_{h},g_{h} \in W_{h}\,.
			\end{equation}
			Therefore we may write
			\begin{equation}
				\langle\tilde{f}_{h}, \tilde{g}_{h}\rangle_{\Omega} = \langle\left(\Phi^{*}_{0}\right)^{-1}\left[f_{h}\right], \left(\Phi^{*}_{0}\right)^{-1}\left[g_{h}\right]\rangle_{\Omega} = \langle f_{h}, g_{h}\rangle_{\Omega_{0}}\,.\label{eq:inner_product_0_curved}
			\end{equation}
			
		\subsubsection{$\tilde{U}_{h}\subset H(\mathrm{div},\Omega)$}
			Any polynomial vector field $\vec{u}_{h}\in U_{h}\subset H(\mathrm{div},\Omega_{0})$, as was seen in Section~\ref{sec:two_dimensional_basis_functions}, \eqref{eq::expansion_edge_polynomials},  can be written as
			\begin{equation}
				\vec{u}_{h} = \sum_{i=1}^{2p(p+1)}u_{i}\,\vec{\epsilon}_{i}^{\,U}\,,
			\end{equation}
			with $\vec{\epsilon}_{i}^{\,U}$, $i=1,\dots,2p(p+1)$, the basis functions of the space $U_{h}$, such that $U_{h}:=\mathrm{span}\{\vec{\epsilon}_{1}^{\,U}(\xi,\eta), \dots, \vec{\epsilon}_{2p(p+1)}^{\,U}(\xi,\eta)\}$, $p$ the polynomial degree, and $(\xi,\eta)\in\Omega_{0}=[-1,1]\times[-1,1]$, as before.
			
			To construct the vector field $\tilde{\vec{u}}_{h}\in\tilde{U}_{h}$ we simply apply the transformation $\left(\Phi^{*}_{1,\mathrm{div}}\right)^{-1}$, \eqref{eq:pullback_inverse_outer_1_forms}
			\begin{equation}
				\tilde{\vec{u}}_{h} :=  \sum_{i=1}^{2p(p+1)}u_{i}\,\left(\Phi^{*}_{1,\mathrm{div}}\right)^{-1}\left[\vec{\epsilon}_{i}^{\,U}\right] = \sum_{i=1}^{2p(p+1)}u_{i}\, \frac{1}{J}\boldsymbol{\mathsf{J}}\left(\vec{\epsilon}_{i}^{\,U}\circ\Phi^{-1}\right) := \sum_{i=1}^{2p(p+1)}u_{i}\,\tilde{\vec{\epsilon}}_{i}^{\,U}\,,
			\end{equation}
			where we have defined the basis function $\tilde{\vec{\epsilon}}_{i}^{\,U}$ of the space $\tilde{U}_{h}$ as
			\begin{equation}
				\tilde{\vec{\epsilon}}_{i}^{\,U} := \frac{1}{J}\boldsymbol{\mathsf{J}}\left(\vec{\epsilon}_{i}^{\,U}\circ\Phi^{-1}\right) \,, \qquad i = 1,\dots, 2p(p+1)\,. \label{eq:basis_functions_1_omega}
			\end{equation}
			The inner product between two vector fields $\tilde{\vec{u}}_{h}, \tilde{\vec{v}}_{h}\in\tilde{U}_{h}$ follows directly
			\begin{equation}
				\langle\tilde{\vec{u}}_{h}, \tilde{\vec{v}}_{h}\rangle_{\Omega} \stackrel{\eqref{eq:basis_functions_1_omega}}{=} \langle\sum_{i=1}^{2p(p+1)}u_{i}\frac{1}{J}\boldsymbol{\mathsf{J}}\left(\vec{\epsilon}_{i}^{\,U}\circ\Phi^{-1}\right), \sum_{j=1}^{2p(p+1)}v_{j}\frac{1}{J}\boldsymbol{\mathsf{J}}\left(\vec{\epsilon}_{j}^{\,U}\circ\Phi^{-1}\right)\rangle_{\Omega} \stackrel{\eqref{eq:continuous_inner_product_curved_geometries}}{=} \sum_{i,j=1}^{2p(p+1)}u_{i}v_{j}\int_{\Omega}\frac{1}{J^{2}}\boldsymbol{\mathsf{J}}\left(\vec{\epsilon}_{i}^{\,U}\circ\Phi^{-1}\right) \cdot \boldsymbol{\mathsf{J}}\left(\vec{\epsilon}_{i}^{\,U}\circ\Phi^{-1}\right)\mathrm{d}\Omega\,.
			\end{equation}
			If we now apply the transformation $\Phi^{*}_{2}$, \eqref{eq:pullback_2_forms}, to the integrand, we transform the integral into an integral in $\Omega_{0}$ without changing its value
			\begin{equation}
				\langle\tilde{\vec{u}}_{h}, \tilde{\vec{v}}_{h}\rangle_{\Omega}  = \sum_{i,j=1}^{2p(p+1)}u_{i}v_{j}\int_{\Omega}\frac{1}{J^{2}}\boldsymbol{\mathsf{J}}\left(\vec{\epsilon}_{i}^{\,U}\circ\Phi^{-1}\right) \cdot \boldsymbol{\mathsf{J}}\left(\vec{\epsilon}_{i}^{\,U}\circ\Phi^{-1}\right)\mathrm{d}\Omega \stackrel{\eqref{eq:pullback_2_forms_integral_invariance}}{=} \sum_{i,j=1}^{2p(p+1)}u_{i}v_{j}\int_{\Omega_{0}}\left(\boldsymbol{\mathsf{J}}\vec{\epsilon}_{i}^{\,U}\right)\cdot\left(\boldsymbol{\mathsf{J}}\vec{\epsilon}_{j}^{\,U}\right)\,\frac{1}{J}\mathrm{d}\Omega_{0}\,.
			\end{equation}
			If we manipulate this last expression we may write it in a more simplified way as
			\begin{equation}\label{eq:inner_product_Uh}
				 \sum_{i,j=1}^{2p(p+1)}u_{i}v_{j}\int_{\Omega_{0}}\left(\boldsymbol{\mathsf{J}}\vec{\epsilon}_{i}^{\,U}\right)\cdot\left(\boldsymbol{\mathsf{J}}\vec{\epsilon}_{j}^{\,U}\right)\,\frac{1}{J}\mathrm{d}\Omega_{0} = \sum_{i,j=1}^{2p(p+1)}u_{i}v_{j}\int_{\Omega_{0}}\left(\vec{\epsilon}_{i}^{\,U}\right)^{\top}\boldsymbol{\mathsf{G}}\,\vec{\epsilon}_{j}^{\,U}\,\mathrm{d}\Omega_{0}\,,
			\end{equation}
			where we have defined
			\begin{equation}
				\boldsymbol{\mathsf{G}} := \frac{1}{J}\,\boldsymbol{\mathsf{J}}^{\top}\boldsymbol{\mathsf{J}}\,.
			\end{equation}
			To further simplify the notation, we introduce an inner product in $U_{h}$, $\langle\cdot,\cdot\rangle : U_{h}\times U_{h} \mapsto \mathbb{R}$ as
			\begin{equation}
				\langle \vec{u}_{h}, \vec{v}_{h}\rangle_{\Omega_{0}} := \int_{\Omega_{0}} \vec{u}_{h}^{\top}\boldsymbol{\mathsf{G}}\,\vec{v}_{h} \,\mathrm{d}\Omega_{0}\,, \quad \text{with}\quad \vec{u}_{h}, \vec{v}_{h} \in U_{h}\,.\label{eq:inner_product_integral_1_curved}
			\end{equation}
			Therefore we may write
			\begin{equation}
				\langle\tilde{\vec{u}}_{h}, \tilde{\vec{v}}_{h}\rangle_{\Omega} = \langle\left(\Phi^{*}_{1,\mathrm{div}}\right)^{-1}\left[\vec{u}_{h}\right], \left(\Phi^{*}_{1,\mathrm{div}}\right)^{-1}\left[\vec{v}_{h}\right]\rangle_{\Omega} = \langle \vec{u}_{h}, \vec{v}_{h}\rangle_{\Omega_{0}}\,.\label{eq:inner_product_1_curved}
			\end{equation}
			
		\subsubsection{$\tilde{Q}_{h}\subset L^{2}(\Omega)$}		
			We have introduced in Section~\ref{sec:two_dimensional_basis_functions} the basis functions for the space $Q_{h}\subset H(\mathrm{rot},\Omega_{0})$ such that
			\[
				Q_{h} := \mathrm{span}\{\epsilon_{1}^{Q}(\xi,\eta), \dots, \epsilon^{Q}_{p^{2}}(\xi,\eta)\}\,,
			\]
			where as before $p$ is the polynomial degree of the basis functions, and $(\xi,\eta)\in\Omega_{0}=[-1,1]\times[-1,1]$. Therefore any function $f_{h}\in Q_{h}$ may be written as a linear combination of the basis functions, \eqref{eq::volume_polynomials_expansion},
			\begin{equation}
				f_{h} = \sum_{i=1}^{p^{2}}f_{i}\,\epsilon_{i}^{Q}\,.
			\end{equation}
			
			If we apply the transformation rule $\left(\Phi^{*}_{1}\right)^{-1}$, \eqref{eq:pullback_inverse_2_forms}, to $f_{h} \in Q_{h}$ we obtain $\tilde{f}_{h}\in \tilde{Q}_{h}$, its polynomial expansion in the physical domain $\Omega$
			\begin{equation}
				\tilde{f}_{h} := \sum_{i=1}^{p^{2}}f_{i}\,\left(\Phi^{*}_{2}\right)^{-1}\left[\epsilon_{i}^{Q}\right] := \sum_{i=1}^{p^{2}}f_{i}\,\left(\epsilon_{i}^{Q}\circ\Phi^{-1}\right)\frac{1}{J} \, := \sum_{i=1}^{p^{2}}f_{i}\,\tilde{\epsilon}_{i}^{Q}\,, \label{eq:2_form_expansion_curved_domain}
			\end{equation}
			where we have defined the basis functions $\tilde{\epsilon}_{i}^{Q}$ of the space $\tilde{Q}_{h}$ as 
			\begin{equation}
				\tilde{\epsilon}_{i}^{Q} := \frac{1}{J}\left(\epsilon_{i}^{Q}\circ\Phi^{-1}\right)\,, \qquad i=1,\dots,p^{2}\,. \label{eq:basis_functions_2_omega}
			\end{equation}
			Therefore, the inner product between two scalar functions $\tilde{f}_{h}, \tilde{g}_{h}\in \tilde{Q}_{h}$ is
			\begin{equation}
				\langle\tilde{f}_{h}, \tilde{g}_{h}\rangle_{\Omega} \stackrel{\eqref{eq:2_form_expansion_curved_domain}}{=} \langle\sum_{i=1}^{p^{2}}f_{i}\,\frac{1}{J}\left(\epsilon_{i}^{Q}\circ\Phi^{-1}\right), \sum_{j=1}^{p^{2}}g_{j}\,\frac{1}{J}\left(\epsilon_{j}^{Q}\circ\Phi^{-1}\right)\rangle_{\Omega} \stackrel{\eqref{eq:continuous_inner_product_curved_geometries}}{=} \sum_{i,j}^{p^{2}}f_{i}\,g_{j}\int_{\Omega}\left(\epsilon_{i}^{Q}\circ\Phi^{-1}\right) \left(\epsilon_{j}^{Q}\circ\Phi^{-1}\right)\,\frac{1}{J^{2}}\mathrm{d}\Omega\,.
			\end{equation}
			As before, we may now apply the transformation $\Phi^{*}_{2}$, \eqref{eq:pullback_2_forms}, to the integrand in order to transform the integral into an integral in $\Omega_{0}$ without changing the value of the integral
			\begin{equation}
				\langle\tilde{f}_{h}, \tilde{g}_{h}\rangle_{\Omega} = \sum_{i,j}^{p^{2}}f_{i}\,g_{j}\int_{\Omega}\left(\epsilon_{i}^{Q}\circ\Phi^{-1}\right) \left(\epsilon_{j}^{Q}\circ\Phi^{-1}\right)\,\frac{1}{J^{2}}\mathrm{d}\Omega \stackrel{\eqref{eq:pullback_2_forms_integral_invariance}}{=} \sum_{i,j}^{p^{2}}f_{i}\,g_{j}\int_{\Omega_{0}}\epsilon_{i}^{Q} \epsilon_{j}^{Q}\,\frac{1}{J}\mathrm{d}\Omega\,.
			\end{equation}
			To simplify the notation, we define an inner product in $Q_{h}$, $\langle\cdot,\cdot\rangle : Q_{h}\times Q_{h} \mapsto \mathbb{R}$ as
			\begin{equation}\label{eq:inner_product_Qh}
				\langle f_{h}, g_{h}\rangle_{\Omega_{0}} := \int_{\Omega_{0}} f_{h}g_{h} \, \frac{1}{J}\mathrm{d}\Omega_{0}\,, \quad \text{with}\quad f_{h},g_{h} \in Q_{h}\,.
			\end{equation}
			Therefore we may write
			\begin{equation}
				\langle\tilde{f}_{h}, \tilde{g}_{h}\rangle_{\Omega} = \langle\left(\Phi^{*}_{2}\right)^{-1}\left[f_{h}\right], \left(\Phi^{*}_{2}\right)^{-1}\left[g_{h}\right]\rangle_{\Omega} = \langle f_{h}, g_{h}\rangle_{\Omega_{0}}\,.\label{eq:inner_product_2_curved}
			\end{equation}

	\subsection{Exact topological relations}
		In Section~\ref{sec:properties_of_the_basis_functions} we showed that the mixed mimetic spectral element basis functions 
based on geometric degrees of freedom enables the exact representation of topological relations. We wish to show here that by using 
these basis functions together with the transformation rules for scalar and vector fields introduced in 
Section~\ref{sec:transformation_rules_scalar_fields} and Section~\ref{sec:transformation_rules_vector_fields}, respectively, we are 
able to exactly represent these topological relations even on non-affine geometries.
		
		To make the demonstration more compact we will use the symbols $\Lambda^{i}$ and $\Lambda_{h}^{i}$, with $i=\{0,1,2\}$, to 
represent the continuum function spaces used in this work and their corresponding discrete counterparts
		\begin{equation}
			\begin{array} {lll}
				\Lambda^{0}(\Omega_{0}) := H(\mathrm{rot},\Omega_{0}), & \Lambda^{1}(\Omega_{0}) := H(\mathrm{div},\Omega_{0}), & \Lambda^{2}(\Omega_{0}) := L^{2}(\Omega_{0}), \\
				\Lambda^{0}_{h}(\Omega_{0}) := W_{h}, & \Lambda^{1}_{h}(\Omega_{0}) := U_{h}, & \Lambda^{2}_{h}(\Omega_{0}) := Q_{h}.
			\end{array}
		\end{equation}
		and
		\begin{equation}
			\begin{array} {lll}
				\Lambda^{0}(\Omega) := H(\mathrm{rot},{\Omega}), & \Lambda^{1}({\Omega}) := H(\mathrm{div},{\Omega}), & \Lambda^{2}({\Omega}) := L^{2}({\Omega}), \\
				\Lambda^{0}_{h}({\Omega}) := \tilde{W}_{h}, & \Lambda^{1}_{h}({\Omega}) := \tilde{U}_{h}, & \Lambda^{2}_{h}({\Omega}) := \tilde{Q}_{h}.
			\end{array}
		\end{equation}
		We will also use the symbol $D^{i}$, with $i=\{0,1\}$, to represent the two differential operators
		\begin{equation}
			\begin{array} {ll}
				D^{0} := \nabla^{\perp}, & D^{1} := \nabla\cdot, \\
				\tilde{D}^{0} := \tilde{\nabla}^{\perp}, & \tilde{D}^{1} := \tilde{\nabla}\cdot.
			\end{array}
		\end{equation}
		
		We wish to demonstrate that even on non-affine meshes, without exact integration, this formulation allows for an exact representation of topological relations of the type
		\begin{equation}
			\tilde{D}^{i}\tilde{f}_{h}^{(i)} = \tilde{g}^{(i+1)}_{h}\,, \label{eq:topological_relation_1}
		\end{equation}
		with $\tilde{f}_{h}^{(i)} \in \Lambda^{i}_{h}(\Omega)$ and $\tilde{g}_{h}^{(i+1)} \in \Lambda^{i+1}_{h}(\Omega)$. We have seen that we can expand the discrete function $\tilde{f}_{h}$ and $\tilde{g}_{h}$ as a linear combination of the basis functions, see \eqref{eq::2d_nodal_polynomials_expansion}, \eqref{eq::expansion_edge_polynomials}, and \eqref{eq::volume_polynomials_expansion}, such that \eqref{eq:topological_relation_1} can be rewritten as
		\begin{equation}
			\tilde{D}^{i}\sum_{n=1}^{\mathrm{dim}\,\Lambda^{i}_{h}(\Omega)}f_{n}\tilde{\epsilon}_{n}^{i} = \sum_{n=1}^{\mathrm{dim}\,\Lambda^{i+1}_{h}(\Omega)}g_{n}\tilde{\epsilon}_{n}^{i+1} \,, \label{eq:topological_relation_2}
		\end{equation}
		where we have removed the tilde from the coefficients $f_{n}$ and $g_{n}$ because these are coefficients associated to the geometric degrees of freedom and, therefore, are invariant under the transformation rules, as shown in Section~\ref{sec:transformation_rules_scalar_fields} and Section~\ref{sec:transformation_rules_vector_fields}. Since we are in a finite element setting, this equation will appear as a weak formulation
		\begin{equation}
			\sum_{n=1}^{\mathrm{dim}\,\Lambda^{i}_{h}(\Omega)}f_{n}\langle\tilde{D}^{i}\tilde{\epsilon}_{n}^{i}, \tilde{\epsilon}_{m}^{i+1}\rangle_{\Omega} = \sum_{n=1}^{\mathrm{dim}\,\Lambda^{i+1}_{h}(\Omega)}g_{n}\langle\tilde{\epsilon}_{n}^{i+1}, \tilde{\epsilon}_{n}^{i+1}\rangle_{\Omega} \,, \quad \mathrm{with}\quad m=\{1,\dots,\mathrm{dim}\,\Lambda^{i+1}_{h}(\Omega)\}\label{eq:topological_relation_3}\,.
		\end{equation}
		We have also seen that the basis functions in $\Omega$ and $\Omega_{0}$ are related by the transformation rules, see \eqref{eq:basis_functions_0_omega}, \eqref{eq:basis_functions_1_omega}, and \eqref{eq:basis_functions_2_omega}, therefore we may rewrite \eqref{eq:topological_relation_3} as
		\begin{multline}
			\sum_{n=1}^{\mathrm{dim}\,\Lambda^{i}_{h}(\Omega)}f_{n}\langle\tilde{D}^{i}\left(\Phi_{i}^{*}\right)^{-1}\left[{\epsilon}_{n}^{i}\right], \left(\Phi_{i+1}^{*}\right)^{-1}\left[{\epsilon}_{m}^{i+1}\right]\rangle_{\Omega} = \sum_{n=1}^{\mathrm{dim}\,\Lambda^{i+1}_{h}(\Omega)}g_{n}\langle\left(\Phi_{i+1}^{*}\right)^{-1}\left[{\epsilon}_{n}^{i+1}\right], \left(\Phi_{i+1}^{*}\right)^{-1}\left[{\epsilon}_{n}^{i+1}\right]\rangle_{\Omega}\,, \\
\quad \mathrm{with}\quad m=\{1,\dots,\mathrm{dim}\,\Lambda^{i+1}_{h}(\Omega_{0})\}\label{eq:topological_relation_4}\,.
		\end{multline}
		Using the commuting relation between the differential operators and the transformation rules, \eqref{eq:commuting_relations_nabla}, yields
		\begin{multline}
			\sum_{n=1}^{\mathrm{dim}\,\Lambda^{i}_{h}(\Omega)}f_{n}\langle\left(\Phi_{i+1}^{*}\right)^{-1}\left[D^{i}{\epsilon}_{n}^{i}\right], \left(\Phi_{i+1}^{*}\right)^{-1}\left[{\epsilon}_{m}^{i+1}\right]\rangle_{\Omega} = \sum_{n=1}^{\mathrm{dim}\,\Lambda^{i+1}_{h}(\Omega)}g_{n}\langle\left(\Phi_{i+1}^{*}\right)^{-1}\left[{\epsilon}_{n}^{i+1}\right], \left(\Phi_{i+1}^{*}\right)^{-1}\left[{\epsilon}_{n}^{i+1}\right]\rangle_{\Omega} \,, \\
\quad \mathrm{with}\quad m=\{1,\dots,\mathrm{dim}\,\Lambda^{i+1}_{h}(\Omega_{0})\}\label{eq:topological_relation_5}\,.
		\end{multline}
		If we now use \eqref{eq:inner_product_0_curved}, \eqref{eq:inner_product_1_curved}, and \eqref{eq:inner_product_2_curved}, we obtain
		\begin{equation}
			\sum_{n=1}^{\mathrm{dim}\,\Lambda^{i}_{h}(\Omega)}f_{n}\langle D^{i}{\epsilon}_{n}^{i}, \epsilon_{m}^{i+1}\rangle_{\Omega_{0}} = \sum_{n=1}^{\mathrm{dim}\,\Lambda^{i+1}_{h}(\Omega)}g_{n}\langle\epsilon_{n}^{i+1}, \epsilon_{n}^{i+1}\rangle_{\Omega_{0}} \,, \quad \mathrm{with}\quad m=\{1,\dots,\mathrm{dim}\,\Lambda^{i+1}_{h}(\Omega_{0})\}\label{eq:topological_relation_6}\,.
		\end{equation}
		These inner products correspond to integrals involving geometric terms, see \eqref{eq:inner_product_Wh}, \eqref{eq:inner_product_integral_1_curved}, and \eqref{eq:inner_product_Qh}, therefore may only be approximated by numerical quadrature. This means that these relations may only be approximated. Hence, their exact topological nature is lost. 
		
		Within the mixed mimetic spectral element formulation we may use the incidence matrices to exactly express $D^{i}{\epsilon}_{n}^{i}$ as a linear combination of the basis ${\epsilon}_{n}^{i+1}$, see \eqref{eq:hilbert_subcomplex_basis_W} and \eqref{eq:hilbert_subcomplex_basis_U}. Using this in \eqref{eq:topological_relation_6} gives
		\begin{equation}
			\sum_{n,j=1}^{\mathrm{dim}\,\Lambda^{i}_{h}(\Omega),\, \mathrm{dim}\,\Lambda^{i+1}_{h}(\Omega)}\mathsf{E}_{jn}^{i+1,i}f_{n}\langle\epsilon_{j}^{i+1}, \epsilon_{m}^{i+1}\rangle_{\Omega_{0}} = \sum_{n=1}^{\mathrm{dim}\,\Lambda^{i+1}_{h}(\Omega)}g_{n}\langle\epsilon_{n}^{i+1}, \epsilon_{n}^{i+1}\rangle_{\Omega_{0}} \,, \quad \mathrm{with}\quad m=\{1,\dots,\mathrm{dim}\,\Lambda^{i+1}_{h}(\Omega_{0})\}\label{eq:topological_relation_7}\,.
		\end{equation}
		In compact matrix notation \eqref{eq:topological_relation_7} becomes
		\begin{equation}
			\boldsymbol{\mathsf{M}}^{i+1}\boldsymbol{\mathsf{E}}^{i+1,i} \boldsymbol{f} = \boldsymbol{\mathsf{M}}^{i+1}\boldsymbol{g}\,,
		\end{equation}
		with $\mathsf{M}^{i+1}_{j,m} := \langle\epsilon_{j}^{i+1}, \epsilon_{m}^{i+1}\rangle_{\Omega_{0}}$. Since $\boldsymbol{\mathsf{M}}^{i+1}$ is an invertible matrix, we can simply eliminate it from both sides of \eqref{eq:topological_relation_7}, yielding once more the exact topological relations
		\begin{equation}
			\boldsymbol{\mathsf{E}}^{i+1,i} \boldsymbol{f} = \boldsymbol{g}\,.
		\end{equation}
		This means that by using the mixed mimetic spectral element method together with the transformation rules discussed before, 
we are able to exactly represent the topological relations on non-affine geometries. Moreover, this may be achieved without exact integration.

\section{Solution of the shallow water equations}\label{sec::sw}

The shallow water equations are given in rotational form for the velocity $\vec u$ and the fluid depth $h$ as
\begin{subequations}\label{sw_cont}
\begin{align}
\frac{\partial\vec u}{\partial t} &= -(\omega + f)\times\vec u - \nabla\Big(\frac{1}{2}\vec u\cdot\vec u + gh\Big) - c_0\triangle^2\vec u\label{sw_cont_mom},\\
\frac{\partial h}{\partial t} &= -\nabla\cdot\vec F, \label{sw_cont_mas}\
\end{align}
\end{subequations}
where $g = 9.80616 \mathrm{m}/\mathrm{s}^{2}$ is the acceleration due to the earth's gravity, $f = 2\Omega\sin(\phi)$ is the Coriolis term 
(for angular frequency $\Omega = 7.292\times 10^{-5}\mathrm{s}^{-1}$ and latitude $\phi$), $c_0$ is the biharmonic viscosity coefficient 
and $\triangle$ is the Laplacian operator. The vorticity $\omega$ and the mass flux $\vec F$ are given respectively by the diagnostic equations

\begin{subequations}\label{diag_cont}
\begin{align}
\omega &= \nabla\times\vec u \label{diag_cont_w},\\
\vec F &= h\vec u. \label{diag_cont_mf}\
\end{align}
\end{subequations}

Our formulation here differs from that in LPG18 and other mimetic constructions of the shallow water equations
\cite{AL81,RTKS10,MC14}. This is because we are not here concerned with the conservation of potential enstrophy 
$hq^2$, where $q$ is the potential vorticity $q = (\omega + f)/h$. In LPG18 we showed that potential enstrophy 
conservation is contingent on exact integration due to the requirement that the chain rule for the relation 
$\nabla q_h^2 = 2q_h\nabla q_h$ hold in the discrete form. Here we abandon this aspiration because (i) it allows 
us to use inexact \secondRev{GLL } quadrature,  leading to diagonal mass matrices for basis functions in $W_h$ and 
(ii) because our geometry is non-affine and accounts for higher order curvature of the sphere, exact integration 
of inner products including metric terms which account for this curvature would require an extremely high 
quadrature order, which would be prohibitively expensive to compute.

\secondRev{We also note that we have not accounted for bottom topography in the above formulation, however the
inclusion of this term via an expansion of trial functions in $Q_h$, as was done in LPG18, is fairly straight forward.}

Before discretizing the system \eqref{sw_cont} \eqref{diag_cont}, we first express these in the continuous weak
form, by multiplying \eqref{sw_cont_mom} by $\vec\nu\in H(\mathrm{div},\Omega)$, \eqref{diag_cont_w} by 
$\eta\in H(\mathrm{rot},\Omega)$, and \eqref{diag_cont_mf} by 
\secondRev{$\vec\mu\in H(\mathrm{div},\Omega)$}, 
and integrating over the
domain $\Omega$, using standard inner product notation of the form $\langle f,g\rangle_{\Omega} = \int_{\Omega}f\cdot g\mathrm{d}\Omega$.
Note that we have not expressed \eqref{sw_cont_mas} in
the weak form, since the divergence theorem is satisfied point wise in the strong form, and that the biharmonic
viscosity has been omitted (this will be addressed in Section \ref{sec::biharm}). This gives

\begin{subequations}\label{sw_wf}
\begin{align}
\Big\langle\vec\nu,\frac{\partial\vec u}{\partial t}\Big\rangle_{\Omega} &= -\langle\vec\nu,(\omega + f)\times\vec u\rangle_{\Omega} + 
\Big\langle\nabla\cdot\vec\nu,\frac{1}{2}\vec u\cdot\vec u + gh\Big\rangle_{\Omega} && \forall\vec\nu\in H(\mathrm{div},\Omega), \label{sw_wf_mom}\\
\frac{\partial h}{\partial t} &= -\nabla\cdot\vec F, \label{sw_wf_mas}\\
\langle\eta,\omega\rangle_{\Omega} &= -\langle\nabla^{\perp}\eta,\vec u\rangle_{\Omega} && \forall\eta\in H(\mathrm{rot},\Omega), \label{sw_wf_w}\\
\langle\vec\mu,\vec F\rangle_{\Omega} &= \langle\vec\mu,h\vec u\rangle_{\Omega} && \forall\vec\mu\in H(\mathrm{div},\Omega), \label{sw_wf_mf}
\end{align}
\end{subequations}
where we have applied the continuous form of the adjoint relations in \eqref{adjt_div_grad} and \eqref{adjt_rot_curl} 
to the differential operators in \eqref{sw_wf_mom} and \eqref{sw_wf_w} respectively.

As in LPG18 we discretize \eqref{sw_wf} with $w_h, f_h\in W_h$, $\vec u_h, \vec F_h\in U_h$
and $h_h\in Q_h$. Note that in LPG18 we also solved for the kinetic energy per unit mass $K_h\in Q_h$, however
here we directly project $\vec u_h\cdot\vec u_h$ onto the test functions $\epsilon_i^{\,Q}$.
Multiplying \eqref{sw_cont_mom} and \eqref{diag_cont_mf} by $\vec\epsilon_i^{\,U}$ and \eqref{diag_cont_w} by 
$\epsilon_i^{\,Q}$, and integrating over the domain $\Omega$ leads to the discrete system

\begin{subequations}{\label{sw_disc}}
\begin{align}
  \sum_{i=1}^{d_{U}}\ip{\vec{\epsilon}^{\,U}_{j}}{\vec{\epsilon}^{\,U}_{i}}\frac{\mathrm{d}u_{i}}{\mathrm{d}t} &=
  - \sum_{i=1}^{d_{U}}\ip{\vec{\epsilon}_j^{\,U}}{(\omega_{h} + f_h)\times\vec{\epsilon}_{i}^{\,U}}u_i
  + g\sum_{i,k=1}^{d_{Q}}(\mathsf{E}^{2,1}_{k,j})^{\top}\ip{\epsilon_{k}^{Q}}{\epsilon_{i}^{Q}}h_i 
  + \frac{1}{2}\sum_{i=1}^{d_{U}}\sum_{k=1}^{d_{Q}}(\mathsf{E}^{2,1}_{k,j})^{\top}\ip{\epsilon_{k}^{Q}}{\vec{u}_{h}\cdot\vec{\epsilon}_{i}^{\,U}}u_i,
\label{mom_cont_weak_discrete_expansion} \\
  \frac{\mathrm{d}h_{j}}{\mathrm{d}t} &= -\sum_{i=1}^{d_{U}}\mathsf{E}^{2,1}_{j,i}F_i,
\label{mas_cont_weak_discrete_expansion} \\
  \sum_{i=1}^{d_{W}}\ip{\epsilon_{j}^{\,W}}{\epsilon_{i}^{\,W}}\omega_i &= 
  - \sum_{i,k=1}^{d_{U}} (\mathsf{E}_{k,j}^{1,0})^{\top}\ip{\vec{\epsilon}_{k}^{\,U}}{\vec{\epsilon}_{i}^{\,U}}u_i +
  \sum_{i=1}^{d_{W}}\ip{\epsilon_{j}^{\,W}}{\epsilon_{i}^{\,W}}f_i,
\label{eq:definition_potential_vorticity_weak_discrete_expansion} \\
  \sum_{i=1}^{d_{U}}\ip{\vec{\epsilon}_{j}^{\,U}}{\vec{\epsilon}_{i}^{\,U}}F_i &=
  \sum_{i=1}^{d_U}\ip{\vec{\epsilon}_{j}^{\,U}}{h_h\vec{\epsilon}_{i}^{\,U}}u_i,
\label{eq:definition_h_flux_weak_discrete_expansion}
\end{align}
\end{subequations}
where the differential operators in \eqref{sw_disc} are derived from the weak form relations \eqref{grad_wf} and \eqref{curl_wf}.
Note that we have not applied a Galerkin projection for \eqref{mas_cont_weak_discrete_expansion}. This is because the divergence 
operator holds point-wise in the strong form as given in \eqref{eq:hilbert_subcomplex_basis_U}. Note also that we have omitted
the biharmonic viscosity term, as this will be discussed later.

We may alternatively express \eqref{sw_disc} in matrix form for each element. First we define
a set of matrices corresponding to the evaluation of basis functions at quadrature points as

\begin{equation}\label{basis_mats}
\bm{A}: \mathsf{A}_{i,j} = \epsilon^{\,W}_j(\vec q_i(\xi_a,\eta_b))\quad
\bm{B}: \mathsf{B}_{i,j} = \vec\epsilon^{\,U}_j(\vec q_i(\xi_a,\eta_b))\quad
\bm{C}: \mathsf{C}_{i,j} = \epsilon^{\,Q}_j(\vec q_i(\xi_a,\eta_b))
\end{equation}
where $\vec q_i(\xi_a,\eta_b)$ is the $i^{\mathrm{th}}$ \secondRev{GLL } quadrature point in the canonical domain of the 
element $\Omega=[-1, 1]\times [-1, 1]\subset \mathbb{R}^{2}$, for $i = b(p + 1) + a + 1$, and $j$ is the basis function 
index. In addition to these we also introduce the matrix $\bm{P}$, which gives the values of the $H(\mathrm{div};\Omega)$ 
form of the Piola transformation \cite{RKL09,RHCM13,NSC16}, as derived in \eqref{eq:pullback_inverse_outer_1_forms}, at 
each quadrature point, and is structured as

\begin{equation}\label{piola_disc}
\bm{P} =
\begin{bmatrix}
\frac{\cos(\phi)\theta,_{\xi}|_{q_1}}{det(J|_{q_1})} & \frac{\cos(\phi)\theta,_{\eta}|_{q_1}}{det(J|_{q_1})} & 0 & 0 & \dots & 0 & 0\\
\frac{\phi,_{\xi}|_{q_1}}{det(J|_{q_1})}             & \frac{\phi,_{\eta}|_{q_1}}{det(J|_{q_1})}             & 0 & 0 & \dots & 0 & 0\\
0 & 0 & \frac{\cos(\phi)\theta,_{\xi}|_{q_2}}{det(J|_{q_2})} & \frac{\cos(\phi)\theta,_{\eta}|_{q_2}}{det(J|_{q_2})} & \dots & 0 & 0\\
0 & 0 & \frac{\phi,_{\xi}|_{q_2}}{det(J|_{q_2})}             & \frac{\phi,_{\eta}|_{q_2}}{det(J|_{q_2})}             & \dots & 0 & 0\\
\vdots                    & \vdots                       & \vdots                      & \vdots                      & \ddots & 0 & 0\\
0 & 0 & 0 & 0 & 0 & \frac{\cos(\phi)\theta,_{\xi}|_{q_{(p+1)^2}}}{det(J|_{q_{(p+1)^2}})} & \frac{\cos(\phi)\theta,_{\eta}|_{q_{(p+1)^2}}}{det(J|_{q_{(p+1)^2}})} \\
0 & 0 & 0 & 0 & 0 & \frac{\phi,_{\xi}|_{q_{(p+1)^2}}}{det(J|_{q_{(p+1)^2}})}             & \frac{\phi,_{\eta}|_{q_{(p+1)^2}}}{det(J|_{q_{(p+1)^2}})}             \\
\end{bmatrix},
\end{equation}
where
\begin{equation}\label{jacobian}
J|_{q_i} = 
\begin{bmatrix}
\cos(\phi)\theta,_{\xi}|_{q_i} & \cos(\phi)\theta,_{\eta}|_{q_i} \\
\phi,_{\xi}|_{q_i}             & \phi,_{\eta}|_{q_i}             \\
\end{bmatrix}
\end{equation}
is the Jacobian evaluated at quadrature point $q_i(\xi,\eta)$, and $det(J|_{q_i})$ is its determinant.
Note that the $\vec e_{\theta}$ terms have been scaled by $\cos(\phi)$, since 
$\mathrm{d}\vec\theta = \cos(\phi)\mathrm{d}\theta\vec e_{\theta} + \mathrm{d}\phi\vec e_{\phi}$
for $\vec\theta = (\theta\vec e_{\theta},\phi\vec e_{\phi})$, where $\vec e_{\theta}$ is the zonal 
direction and $\vec e_{\phi}$ is the meridional direction.
The specifics of the implementation of the Jacobian terms can be found in \cite{GTUOL14}.

We also introduce a diagonal matrix $\bm{D}$ which contains the Jacobian determinant evaluated at each quadrature point $q_i$,
and a second diagonal matrix $\bm{W}$, which corresponds to the weights of the quadrature points within the canonical element.
Using these operators, local degrees of freedom 
$\boldsymbol{\omega}:=[\omega_{1},\dots,\omega_{d_{W}}]^{\top}$.
$\boldsymbol{u}:=[u_{1},\dots,u_{d_{U}}]^{\top}$, and
$\boldsymbol{h}:=[h_{1},\dots,h_{d_{Q}}]^{\top}$,
may be mapped to global spherical coordinates via \eqref{eq:pullback_inverse_0_forms}, \eqref{eq:pullback_inverse_outer_1_forms} 
and \eqref{eq:pullback_inverse_2_forms} respectively as

\begin{subequations}
\begin{align}
\boldsymbol{\omega}_g &= \bm{A}\boldsymbol{\omega},\\
\boldsymbol{u}_g &= \bm{P}\bm{B}\boldsymbol{u} = \hat{\bm{B}}\boldsymbol{u},\\
\boldsymbol{h}_g &= \bm{D}^{-1}\bm{C}\boldsymbol{h} = \hat{\bm{C}}\boldsymbol{h}.
\end{align}
\end{subequations}
Note that no mapping from local to global coordinates is required for variables in $W_h$.

For $\bm{P}$ in \eqref{piola_disc} the vector component represents the minor index for each row and column and the quadrature
point is the major index. For all other matrices operating on vector fields in $U_h$ this convention also applies, that is the 
minor index is the vector component and the major index corresponds to the quadrature point or basis function. Having established 
this convention we introduce four other pieces of notation: a subscript $2$, which when applied to a diagonal matrix as $\bm{X}_2$
creates a minor index of length 2 by repeating the entries of the scalar field; a subscript $d$, which when
applied to a vector as $\bm{x}_d$ converts this to a diagonal matrix; a superscript $\perp$, which when applied to a matrix as 
$\bm{X}^{\perp}$ or a vector as $\bm{x}^{\perp}$ rotates the minor index components by $\pi/2$ radians by the rotational matrix

\begin{equation}
\begin{bmatrix}
0 & -1 \\
1 &  0 \\
\end{bmatrix};
\end{equation}
and the operator $\bar\cdot$, which takes an inner product over the minor indices.
Having outlined this notation, we can represent \eqref{sw_disc} in matrix form including the metric 
transformations as

\begin{subequations}\label{sw_mat}
\begin{align}
\hat{\bm{B}}^{\top}\bm{D}_2\bm{W}_2\hat{\bm{B}}\frac{\partial\boldsymbol{u}}{\partial t} &= 
-\hat{\bm{B}}^{\top}\bm{D}_2\bm{W}_2(\bm{A}(\boldsymbol{\omega} + \boldsymbol{f}))_{d,2}\hat{\bm{B}}^{\perp}\boldsymbol{u}^{\perp} + 
(\bm{E}^{2,1})^{\top}\hat{\bm{C}}^{\top}\bm{D}\bm{W}(0.5(\hat{\bm{B}}\boldsymbol{u})_d\bar\cdot\hat{\bm{B}}\boldsymbol{u} + g\hat{\bm{C}}\boldsymbol{h}),\label{sw_mat_mom}\\
\frac{\partial\boldsymbol{h}}{\partial t} &= -\bm{E}^{2,1}\boldsymbol{F}, \label{sw_mat_mas}\\
\bm{A}^{\top}\bm{D}\bm{W}\bm{A}\boldsymbol{\omega} &= -(\bm{E}^{1,0})^{\top}\hat{\bm{B}}^{\top}\bm{D}_2\bm{W}_2\hat{\bm{B}}\boldsymbol{u},\label{diag_mat_w}\\
\hat{\bm{B}}^{\top}\bm{D}_2\bm{W}_2\hat{\bm{B}}\boldsymbol{F} &= \hat{\bm{B}}\bm{D}_2\bm{W}_2(\hat{\bm{C}}\boldsymbol{h})_{d,2}\hat{\bm{B}}\boldsymbol{u},\label{diag_mat_mf}
\end{align}
\end{subequations}
where $\boldsymbol{F}:=[F_{1},\dots,F_{d_{U}}]^{\top}$, and the mass matrices on the left hand sides of \eqref{sw_mat_mom},
and \eqref{diag_mat_mf} are the matrix forms of the inner product relation \eqref{eq:inner_product_Uh}, and the left hand
side of \eqref{diag_mat_w} is the matrix form of \eqref{eq:inner_product_Wh}.
Note that since we are using inexact quadrature, such that the quadrature points are collocated with the nodal points, 
the matrix $\bm{A}$ is diagonal, and so taking the transpose of this in \eqref{diag_mat_w} is unnecessary, however
we keep this in for completeness.

\subsection{Conservation}

In LPG18 we followed the work of \cite{AL81,MC14} to demonstrate the conservation properties of the discrete form 
of the shallow water equations using mixed mimetic spectral elements. Here we simply re-state the principle criteria 
for conservation of the various moments.

\begin{itemize}
\item Mass, $\frac{\mathrm{d}}{\mathrm{d}t}\int h_h\mathrm{d}\Omega = 0$: 
Point-wise conservation due to the strong form representation of the divergence theorem \eqref{eq:hilbert_subcomplex_basis_U}.

\item Vorticity, $\frac{\mathrm{d}}{\mathrm{d}t}\int w_h\mathrm{d}\Omega = 0$: 
Weak form conservation due to the annihilation of the gradient by the curl \eqref{curl_grad_wf},
and the mapping of basis functions in $W_h$ into basis functions in $U_h$ by the discrete \emph{rot} operator \eqref{eq:hilbert_subcomplex_basis_W}.

\item Energy, $\frac{\mathrm{d}}{\mathrm{d}t}\ip{h_h}{\frac{1}{2}\vec u_h\cdot\vec u_h + \frac{1}{2}gh_h} = 0$:
Orthogonality of $\vec u_h$ and $(\omega_h + f_h)\times\vec u_h$
(or alternatively the anti-symmetry of the convective operator in \eqref{sw_mat_mom} 
\cite{ER17}) and the adjoint relationship between \emph{grad} and \emph{div} \eqref{adjt_div_grad} \cite{Salmon04}.
\end{itemize}

\subsection{Biharmonic viscosity}\label{sec::biharm}

We also apply a biharmonic viscosity operator to \eqref{sw_mat_mom} of the form $c_0\triangle_h^2\vec u_h$, where
$c_0 = 0.0718\Delta x^{3.2}$, as derived from previous spectral element simulations \cite{GTUOL14}, and $\Delta x$
is the average nodal grid spacing on the cubed sphere. The Laplacian operator is applied via the vector identity
$\triangle\vec u = \nabla(\nabla\cdot\vec u) + \nabla^{\perp}(\nabla\times\vec u)$.
The continuous weak form is given via the Helmholtz decomposition \cite{KPG11} as 

\begin{subequations}
\begin{align}
\langle\eta,\omega\rangle_{\Omega} &= -\langle\nabla^{\perp}\eta,\vec u\rangle_{\Omega} && \forall\eta\in H(\mathrm{rot},\Omega),\label{biharm_1}\\
\vec r &= \nabla^{\perp}\omega,\label{biharm_2}\\
\delta &= \nabla\cdot\vec u,\label{biharm_3}\\
\langle\vec\nu,\vec d\rangle_{\Omega} &= -\langle\nabla\cdot\vec\nu,\delta\rangle_{\Omega}&& \forall\vec\nu\in H(\mathrm{div},\Omega),\label{biharm_4}\\
\triangle\vec u &= \vec r + \vec d.\label{biharm_5}
\end{align}
\end{subequations}
Note that \eqref{biharm_2} and \eqref{biharm_3} are expressed in the strong form since the \emph{rot} and \emph{div}
operators are applied point wise.
The equivalent discrete matrix form of the Laplacian is then given as

%We also apply a biharmonic viscosity operator to \eqref{sw_mat_mom} of the form $c_0\triangle_h^2\vec u_h$, where
%$c_0 = 0.0718\Delta x^{3.2}$, as derived from previous spectral element simulations \cite{GTUOL14}, and $\Delta x$
%is the average nodal grid spacing on the cubed sphere.
%The application of the Laplacian operator $\triangle\vec u_h$ is achieved using the discrete Helmholtz decomposition 
%\cite{KPG11} of $\vec u_h$ as

\begin{subequations}\label{biharmonic_visc_disc}
\begin{align}
\bm{A}^{\top}\bm{D}\bm{W}\bm{A}\boldsymbol{\omega} &= -(\bm{E}^{1,0})^{\top}\hat{\bm{B}}^{\top}\bm{D}_2\bm{W}_2\hat{\bm{B}}\boldsymbol{u},\\
\boldsymbol{r} &= \bm{E}^{1,0}\boldsymbol{\omega},\\
\boldsymbol{\delta} &= \bm{E}^{2,1}\boldsymbol{u},\\
\hat{\bm{B}}^{\top}\bm{D}_2\bm{W}_2\hat{\bm{B}}\boldsymbol{d} &= -(\bm{E}^{2,1})^{\top}\hat{\bm{C}}^{\top}\bm{D}\bm{W}\hat{\bm{C}}\boldsymbol{\delta},\\
\triangle_h\boldsymbol{u} &= \boldsymbol{d} + \boldsymbol{r},
\end{align}
\end{subequations}
where $\boldsymbol{d}:=[d_1,\dots,d_{d_U}]^{\top}$, $\boldsymbol{r}:=[r_1,\dots,r_{d_U}]^{\top}$,
and $\boldsymbol{\delta}:=[\delta_1,\dots,\delta_{d_Q}]^{\top}$ is the divergence.
The biharmonic viscosity is then determined by applying this process a second time. 
%The system \eqref{biharmonic_visc_disc} 
%is the discrete form of the vector identity $\triangle\vec u = \nabla(\nabla\cdot\vec u) + \nabla^{\perp}(\nabla\times\vec u)$.

\subsection{Implementation}

We solve the explicit system \eqref{sw_disc} using a stiffly stable second order Runge-Kutta time integrator 
\cite{SO88} of the form

\[y' = y^n - \Delta tG(y^n),\qquad y^{n+1} = y^n - \frac{1}{2}\Delta t(G(y^n) + G(y'))\]
for $y = [\vec u, h]^{\top}$, $G = [(\omega + f)\times\vec u + \nabla(\frac{1}{2}\vec u\cdot\vec u + gh) + c_0\triangle^2\vec u,\nabla\cdot\vec F]^{\top}$.

\secondRev{We generate the } geometry of the element corners on the cubed sphere \secondRev{using } the gnomonic projection in \cite{GPF03}, 
\secondRev{while } the geometry of the internal and edge nodes of the elements \secondRev{are } determined from the Jacobian mapping \cite{GTUOL14}. 
The topological ordering of the faces is the same as for the C-CAM model, a C-grid finite difference primitive equation model on the
cubed sphere \cite{McGregor05}. This ensures that the orientation of the edges is the same for all faces,
and the sign of the Jacobian determinant does not change. 
\secondRev{We use a } simple parallel
decomposition whereby each face of the cubed sphere is divided into $n$ regions in each dimension
such that the total number of processors for the cubed sphere is $6\times n^2$.
\secondRev{The mesh topology and parallel decomposition are generated in a stand alone initialization script and written
to separate files for each processor that are then imported by the main code at runtime. }

The PETSc library \cite{petsc-web-page,petsc-user-ref,petsc-efficient} is used \secondRev{for parallel message passing, 
via the scattering of } data from global to local 
vectors and vice versa, as well as to assemble matrices in parallel and to solve linear systems. 
\secondRev{Equations \eqref{sw_mat_mom} and \eqref{diag_mat_mf} are solved using the GMRES iterative method with a block Jacobi 
preconditioner, with a block size of $2p(p+1)$. Using this solution strategy, these linear systems converge to
a relative tolerance of $10.0^{-16}$ in approximately 
20-30 iterations. We have not extensively explored different solver/preconditioner strategies and use these 
settings purely on account of their generality and robust performance.}

\section{Results}\label{sec::res}

We validate the model using standard test cases \cite{WDHJS92,GSP04}. These are used to validate both the spectral
convergence of errors for analytical solutions, as well as the mimetic conservation properties of the scheme.

\subsection{Williamson test case 2 (exact steady state solution)}

The first test case \cite{WDHJS92} is an exact steady state solution for the inviscid ($c_0 = 0$) rotating
shallow water equations (test case 2). We use this to validate the spectral convergence of errors in the model. 
We run the test for 5 days of simulation time with $\alpha = \pi/4$, such that the flow contains a significant 
meander, in order to avoid the fortuitous cancellation of errors that may occur in a flow that is more closely 
alligned with the mesh ($\alpha=0.0$).

\begin{figure}[!hbtp]
\begin{center}
\includegraphics[width=0.48\textwidth,height=0.36\textwidth]{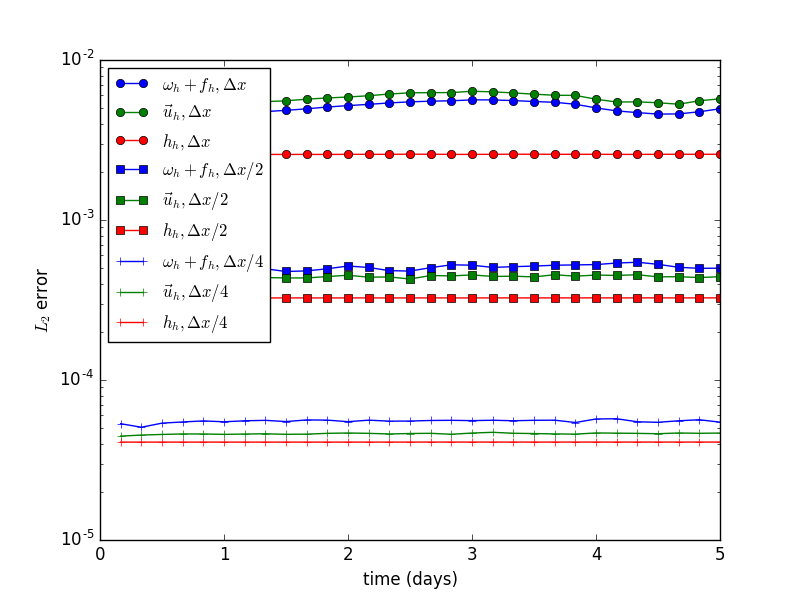}
\includegraphics[width=0.48\textwidth,height=0.36\textwidth]{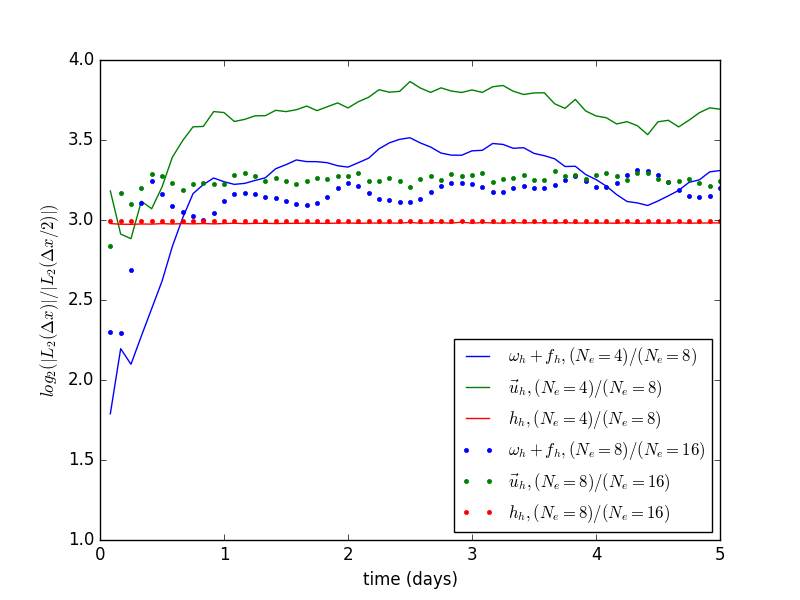}\\
\caption{$L_2$ errors for vorticity, velocity and depth with time for Williamson test case 2 with $p = 3$ (left), 
and $log_2$ of the ratio of $L_2$ errors between $N_e = 4$ $(\Delta t = 240\mathrm{s})$, $N_e = 8$ $(\Delta t = 120\mathrm{s})$
and $N_e = 16$ $(\Delta t = 60\mathrm{s})$ (right). 
Error convergence is approximately $3^{\mathrm{rd}}$ order for all variables.}
\label{fig::williamson_2_1}
\end{center}
\end{figure}

Figure \ref{fig::williamson_2_1} shows both the $L_2$ errors for the absolute vorticity, $\omega_h + f_h$, $\vec u_h$ and $h_h$ (left),
and the convergence of errors between resolutions of $N_e = 4$,  $N_e = 8$ and $N_e = 16$ with time (right), where $N_e$
is the number of elements in each dimension for each face of the cubed sphere. The convergence is determined 
by taking the $log_2$ of the ratio of the errors between higher and lower resolutions, in order to show the 
rate at which these reduce as the resolution is halved.
The error convergence for $p=3$ basis functions are almost precisely $3^{\mathrm{rd}}$ order between nodal 
resolutions of $\Delta x\approx 834\mathrm{km}$, $\Delta x\approx 417\mathrm{km}$ and $\Delta x\approx 207\mathrm{km}$ 
at the equator (the anticipated spectral convergence 
rate for $p = 3$) for $h_h$, while errors for $\omega_h + f_h$ and $\vec u_h$ vary somewhat with time. 
The fact that the error convergence for $\omega_h + f_h$ and $\vec u_h$ is greater than $3^{\mathrm{rd}}$ order after
the initial spin up stage is thought to be because the nodal basis functions $l_i(\xi)$ are a degree higher than 
the edge functions $e_i(\xi)$, and the basis functions for $\vec u_h \in U_h$ are comprised of nodal functions
in one dimension, while the basis functions for $\omega_h + f_h\in W_h$ are composed of nodal basis functions in both 
dimensions. The fact that errors for $\omega_h + f_h$ are somewhat lower than those for $\vec u_h$ despite the use of
nodal basis functions in both dimensions is perhaps because we do not explicitly solve for $\omega_h + f_h$, but rather 
diagnose it in the weak form.

\begin{figure}[!hbtp]
\begin{center}
\includegraphics[width=0.48\textwidth,height=0.36\textwidth]{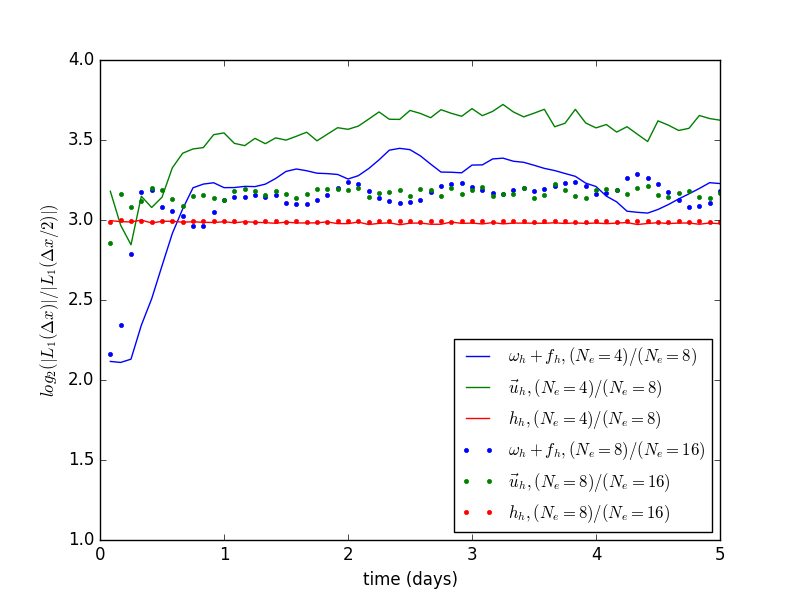}
\includegraphics[width=0.48\textwidth,height=0.36\textwidth]{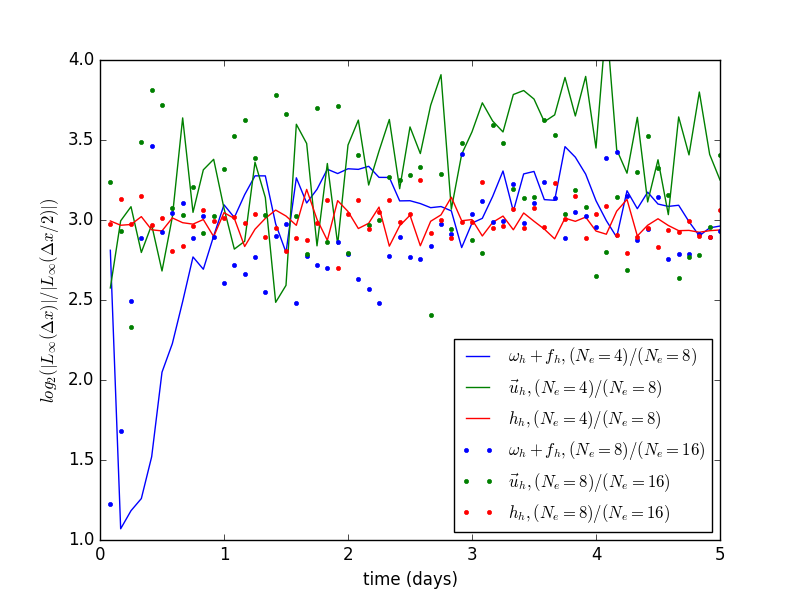}\\
\caption{$log_2$ of the ratio of $p = 3$ errors between $N_e = 4$ $(\Delta t = 240\mathrm{s})$, $N_e = 8$ $(\Delta t = 120\mathrm{s})$
and $N_e = 16$ $(\Delta t = 60\mathrm{s})$ in the $L_1$ norm (left) and the $L_{\infty}$ norm (right) for Williamson test case 2.}
\label{fig::williamson_2_2}
\end{center}
\end{figure}

\begin{figure}[!hbtp]
\begin{center}
\includegraphics[width=0.48\textwidth,height=0.36\textwidth]{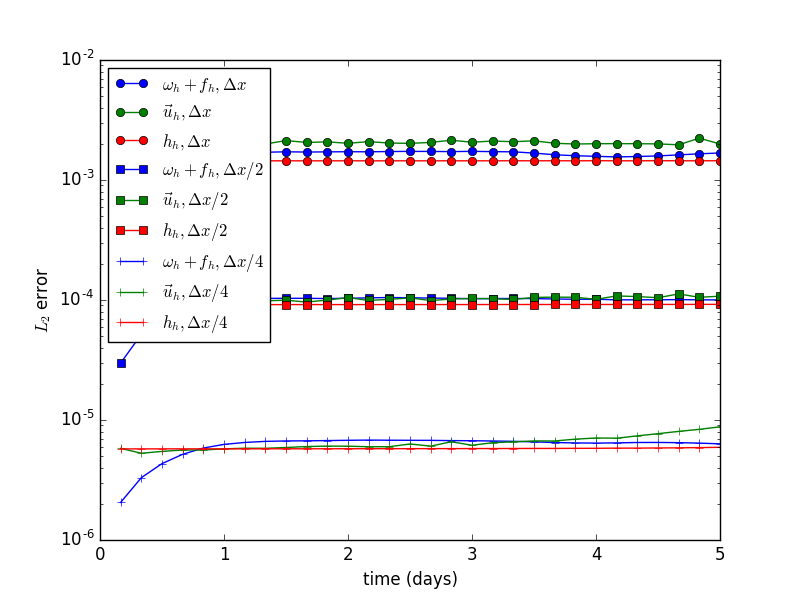}
\includegraphics[width=0.48\textwidth,height=0.36\textwidth]{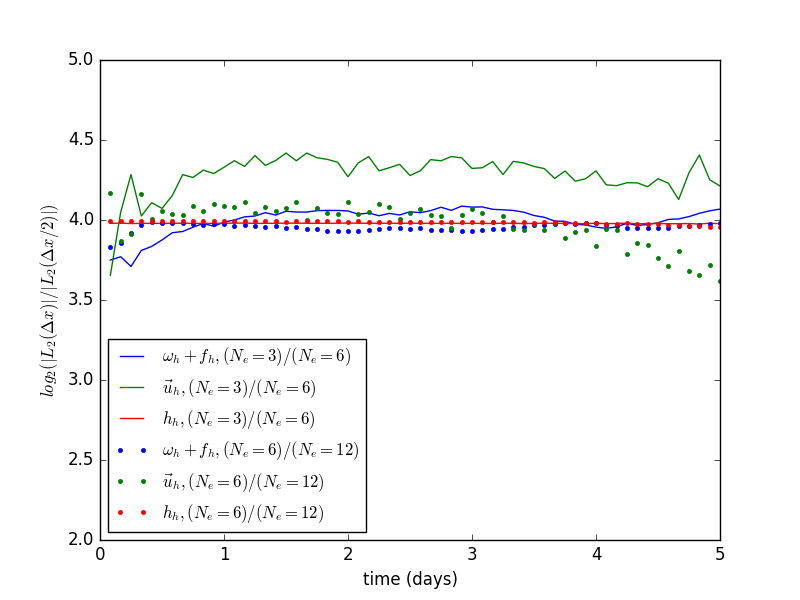}\\
\caption{As for Fig. \ref{fig::williamson_2_1}, however with $p = 4$ and comparison between $N_e = 3$ $(\Delta t = 240\mathrm{s})$,
$N_e = 6$ $(\Delta t = 120\mathrm{s})$ and $N_e = 12$ $(\Delta t = 60\mathrm{s})$.}
\label{fig::williamson_2_3}
\end{center}
\end{figure}

\begin{figure}[!hbtp]
\begin{center}
\includegraphics[width=0.48\textwidth,height=0.36\textwidth]{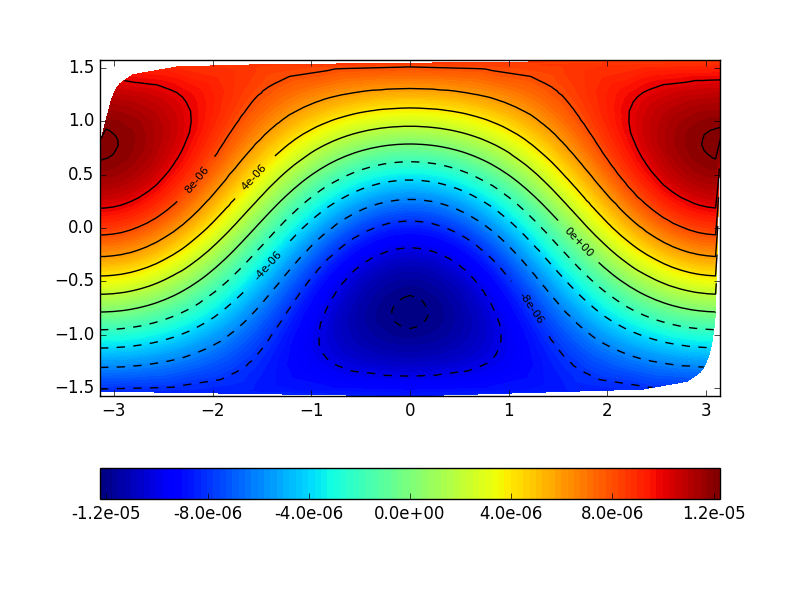}
\includegraphics[width=0.48\textwidth,height=0.36\textwidth]{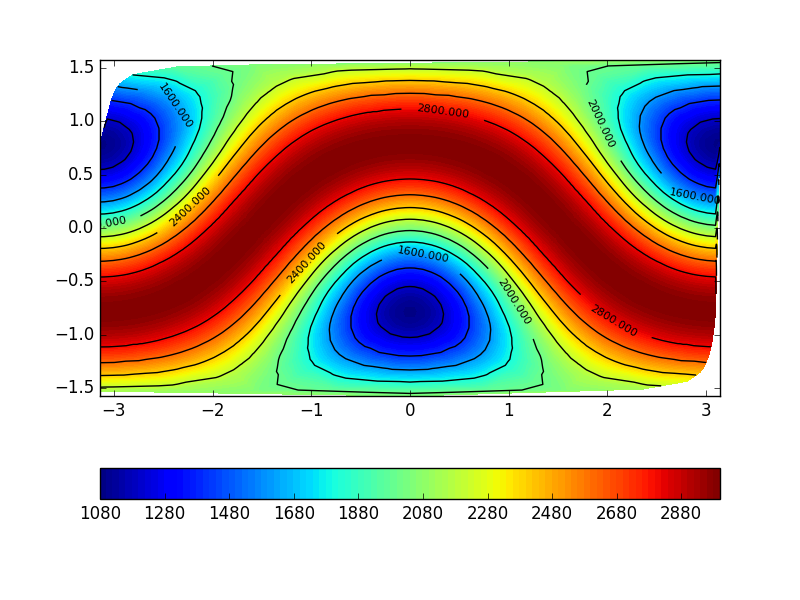}\\
\caption{Vorticity, $\omega_h$ (left), and fluid depth $h_h$ (right), for Williamson test case 2 
at t = 5 days with $p=4$, $N_e=6$, $\Delta t = 120\mathrm{s}$.}
\label{fig::williamson_2_4}
\end{center}
\end{figure}

Results are broadly similar in the $L_1$ norm, with the convergence is somewhat less even in the $L_{\infty}$ norm,
as shown in Fig. \ref{fig::williamson_2_2}. This is to be expected since the cubed sphere mesh is non-uniform, 
so the spread of errors may vary spatially.
The convergence of errors in the $L_2$ norm for $p=4$ between $N_e = 3$, $N_e = 6$ and $N_e = 12$ is shown in Fig. \ref{fig::williamson_2_3}, 
in order to demonstrate that optimal convergence is also preserved at this order. Contour plots for the
vorticity, $\omega_h$ and the fluid depth, $h_h$ with $p = 4$, $N_e = 6$ at day 5 are also shown in Fig. \ref{fig::williamson_2_4}.

\begin{figure}[!hbtp]
\begin{center}
\includegraphics[width=0.48\textwidth,height=0.36\textwidth]{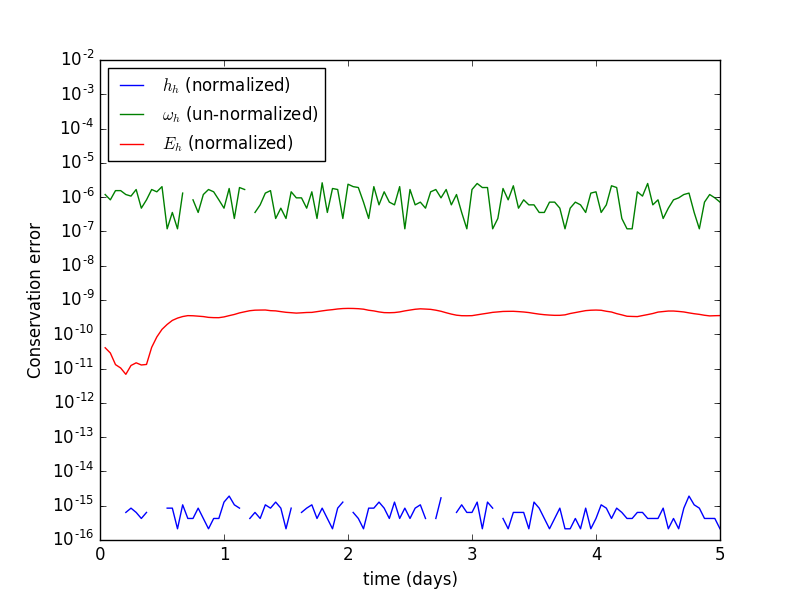}
\includegraphics[width=0.48\textwidth,height=0.36\textwidth]{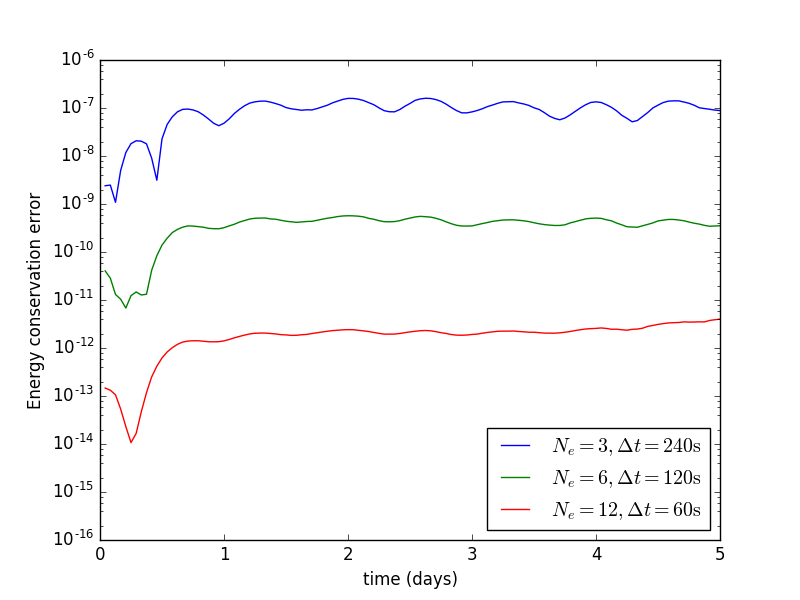}
\caption{Conservation of mass, vorticity and energy for Williamson test case 2 with $p=4$, $N_e = 6$, $\Delta t = 120\mathrm{s}$ (left).
Mass conservation is exact (machine precision), while the (un-normalized) vorticity and energy losses remain bounded but 
are not exact due to the weak representation of vorticity and energy conservation and the use of an 
iterative solver for the momentum equation. Energy conservation will eventually degrade due to
nonlinear cascades to subgrid scales. The convergence of energy conservation errors is also shown for
increasing spatial and temporal resolution (right).}
\label{fig::williamson_2_conservation}
\end{center}
\end{figure}

We also present the mass, vorticity and energy conservation errors in Fig. \ref{fig::williamson_2_conservation}
for test case 2 with $p = 4$, $N_e = 6$, $\Delta t = 120\mathrm{s}$ and $c_0 = 0.0$. Mass conservation is
shown to be exact (machine precision), due to the point-wise preservation of the divergence theorem.
Vorticity conservation errors meanwhile remain bounded and of $\mathcal{O}(10^{-6})$. Note that unlike
the mass and energy conservation errors, vorticity errors are un-normalized, since the global integral of the 
initial vorticity is 0. The fact that vorticity conservation is not preserved to machine precision is most likely
on account of the fact that this is satisfied in the weak form only, and the accuracy of weak form solutions
are bounded by the accuracy of the iterative solver for the $U_h$ mass matrix. The same is true for the
energy conservation errors, which are also bounded but not exact. Note also that with sufficient time 
energy conservation will break down due to nonlinear cascades.

In Fig. \ref{fig::williamson_2_conservation} we also show the convergence of energy conservation errors
with increasing spatial and temporal resolution for $p = 4$. Using a second order Runge-Kutta time integrator
the errors decay at approximately $8^{\mathrm{th}}$ order. The growth of errors due to the nonlinear cascade
grid scales is apparent for the $N_e=12$ case at day 5. While the convergence of energy conservation errors
with temporal and spatial resolution is demonstrated here, we have not shown the convergence of energy 
conservation errors with time step only, since a repeated doubling of the time step required to show this
convergence leads to a time step for which the solution becomes unstable for short simulation times
due to nonlinear cascades in the absence of viscosity. However the fact that the energy conservation errors 
converge as the product of the temporal and spatial orders of the scheme suggests that the convergence of 
these errors with time is behaving correctly.

\begin{figure}[!hbtp]
\begin{center}
\includegraphics[width=0.48\textwidth,height=0.36\textwidth]{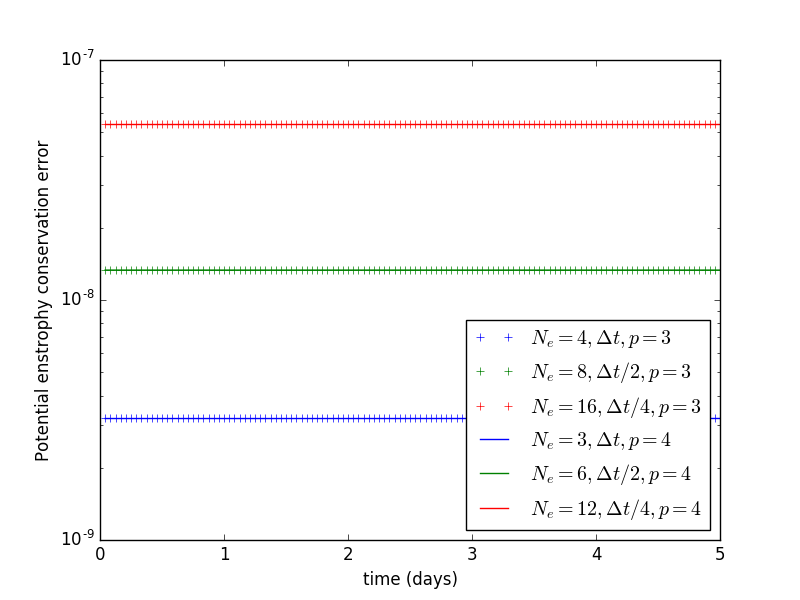}
\includegraphics[width=0.48\textwidth,height=0.36\textwidth]{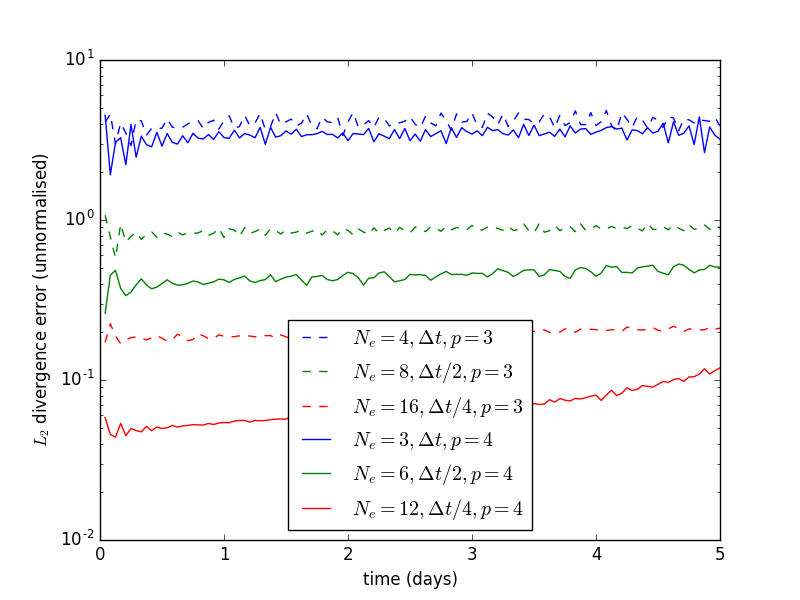}
\caption{Error convergence for potential enstrophy conservation (left) and $L_2$ norm of the divergence 
$\mathsf{E}^{2,1}\boldsymbol u$ (right) for the Williamson test case 2 with $\Delta t = 240s$, 
$\Delta t = 120s$ and $\Delta t = 60s$.}
\label{fig::williamson_2_pe_div}
\end{center}
\end{figure}

\firstRev{We also show the potential enstrophy conservation and the unnormalised divergence
errors in Fig. \ref{fig::williamson_2_pe_div}. As reported in LPG18, potential enstrophy conservation is 
dependent on the preservation of the product rule for quadratic nonlinearities in the discrete form via 
exact quadrature. With inexact GLL quadrature this property is violated, and the potential enstrophy 
conservation errors converge in time only for both $p=3$ and $p=4$ basis functions, such that using a second 
order time integrator these errors reduce by a factor of 4 with a halving of the time step. This is in
contrast to the energy conservation errors shown in Fig. \ref{fig::williamson_2_conservation}, where
we observe energy conservation errors to decrease as the product of the spatial and temporal orders of
the scheme. 

The $L_2$ errors for the divergence, $\mathsf{E}^{2,1}\boldsymbol u$, are shown to converge at one degree 
lower than the polynomial order of the basis functions (second order for $p = 3$ and third order for $p = 4$). 
This result warrants further investigation, as formally the divergence should converge at the same rate
as the polynomial degree in the $L_2$ norm since it is defined on the function space of $Q_h$. The reason for 
the observed convergence rate of the divergence is not currently well understood. This is particularly curious
since the fluid depth, $h_h\in Q_h$ does converge at its anticipated rate as shown in Figs.  
\ref{fig::williamson_2_1},\ref{fig::williamson_2_2} and \ref{fig::williamson_2_3}.}

\subsection{Williamson test case 6 (Rossby-Haurwitz wave)}

The second test case \cite{WDHJS92} is for the evolution of a Rossby-Haurwitz wave (test case 6). Note that
this is an analytical solution of the barotropic vorticity equation, and not the shallow water equations,
and so does not account for the presence of gravity waves. The solution has an east to west group velocity of
$c_g = r_e(R(3+R)\Omega' - 2\Omega)/((1+R)(2+R))\vec e_{\theta}\mathrm{m/s}$, where $r_e = 6.37122\times 10^6\mathrm{m}$ 
is the earth's radius, $R = 4$ is the zonal wave number and $\Omega' = 7.848\times 10^{-6}\mathrm{s}^{-1}$.

\begin{figure}[!hbtp]
\begin{center}
\includegraphics[width=0.48\textwidth,height=0.36\textwidth]{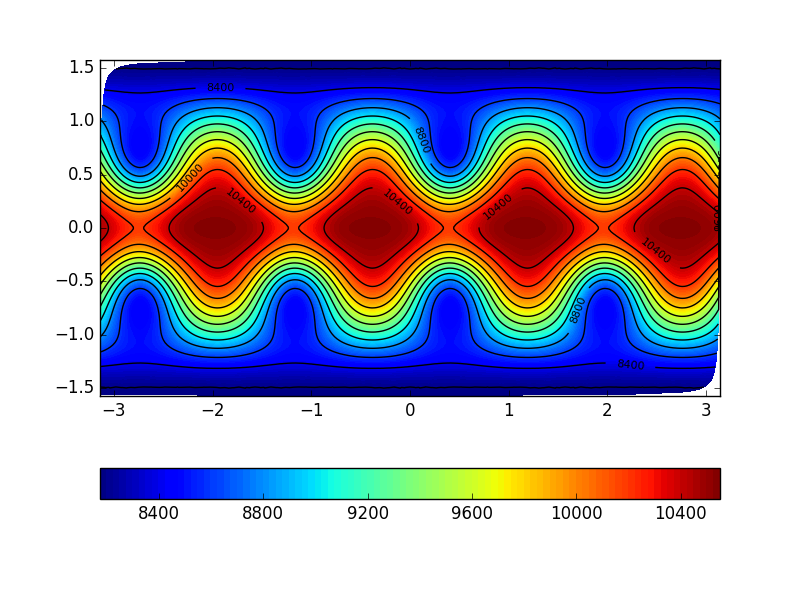}
\includegraphics[width=0.48\textwidth,height=0.36\textwidth]{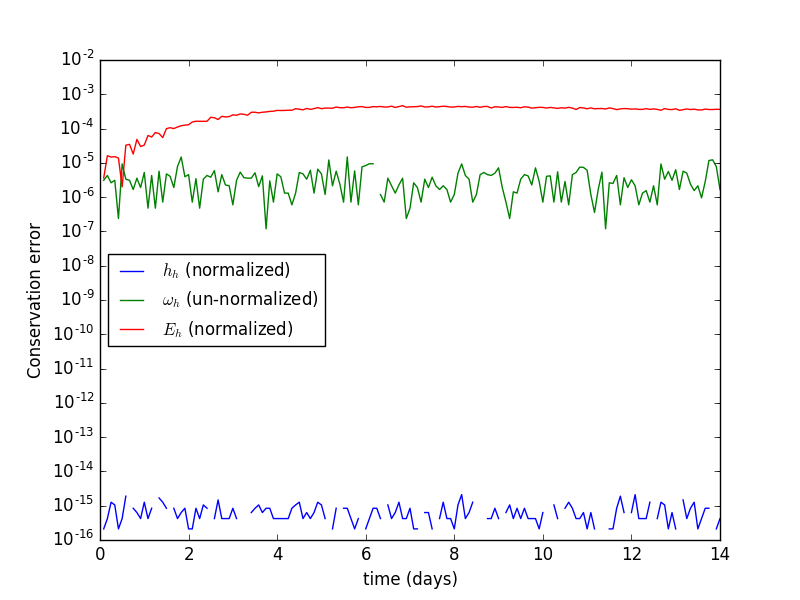}
\caption{Fluid depth for Rossby-Haurwitz wave after 14 days with $p = 3$
$N_e = 32$, $\Delta t = 40\mathrm{s}$ (left).
Conservation of mass, vorticity and energy (right).
Mass conservation is exact (machine precision), while the (un-normalized) vorticity losses remain bounded 
but are not exact
due to the weak representation of vorticity conservation and the use of an iterative solver for the 
momentum equation. Energy losses occur due to the application of biharmonic viscosity.}
\label{fig::williamson_6}
\end{center}
\end{figure}

As for test case 2, we also show the conservation errors for test case 6 in Fig. \ref{fig::williamson_6}.
The mass and vorticity conservation errors are similar to those presented for test case 2 in Fig. 
\ref{fig::williamson_2_conservation}, while energy is not conserved due to the application of 
biharmonic viscosity. With $p=3$ and $N_e=32$, such there are 32 elements in each dimension for each face
of the cubed sphere (corresponding to an approximate resolution of $1^{\circ}$ at the equator), the group
velocity magnitude, $|c_g|$ is approximately $94\%$ if its analytical value. Contours of the depth field $h_h$ after
14 days are consistent with previous A-grid mimetic spectral element results at this resolution \cite{TF10}.

\secondRev{Unlike other compatible schemes with semi-implicit time integration \cite{PG17}, we observe an 
exponentially unstable growth in the solution in the absence of viscosity. 
In LPG18 it was found that with inexact quadrature the product rule can not be preserved for quadratic 
nonlinearities, resulting in a loss of potential enstrophy conservation, leading to an exponential 
growth in the forward cascade of potential enstrophy at grid scales in the absence of viscosity. 
We therefore use biharmonic viscosity not just as
a subgrid turbulence scheme, but also to suppress model instability.}

\subsection{Galewsky test case (barotropic instability of a mid-latitude jet)}

\begin{figure}[!hbtp]
\begin{center}
\includegraphics[width=0.90\textwidth,height=0.45\textwidth]{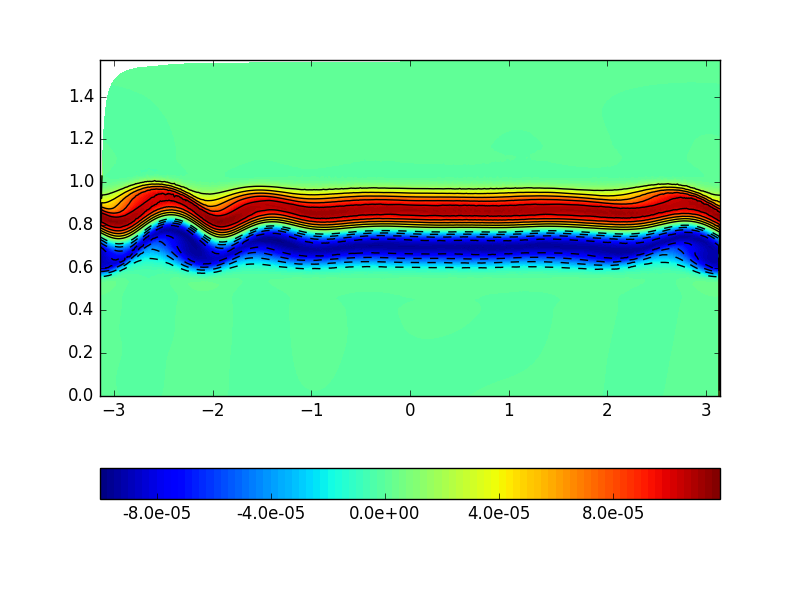}
\caption{Vorticity, $\omega_h$ for the Galewsky test case at day 4 with $p=3$, $N_e=32$ and
$\Delta t = 40\mathrm{s}$. The peak of the instability and contours match well with \cite{GSP04}. 
Contours separation is $2.0\times 10^{-5}\mathrm{s}^{-1}$, with the $\omega_h=0$ contour omitted. 
Only the northern hemisphere is shown.}
\label{fig::galewsky_1}
\end{center}
\end{figure}

\begin{figure}[!hbtp]
\begin{center}
\includegraphics[width=0.90\textwidth,height=0.45\textwidth]{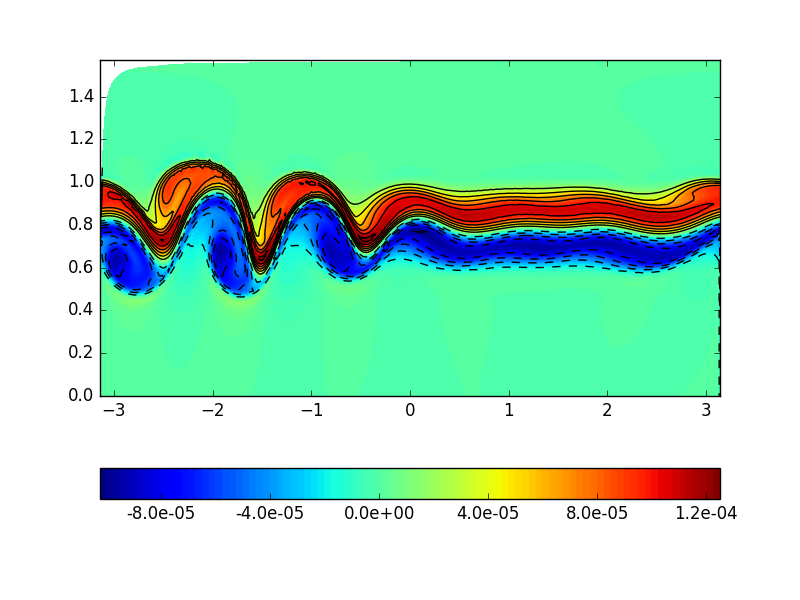}
\caption{As for Fig. \ref{fig::galewsky_1}, day 5.}
\label{fig::galewsky_2}
\end{center}
\end{figure}

\begin{figure}[!hbtp]
\begin{center}
\includegraphics[width=0.90\textwidth,height=0.45\textwidth]{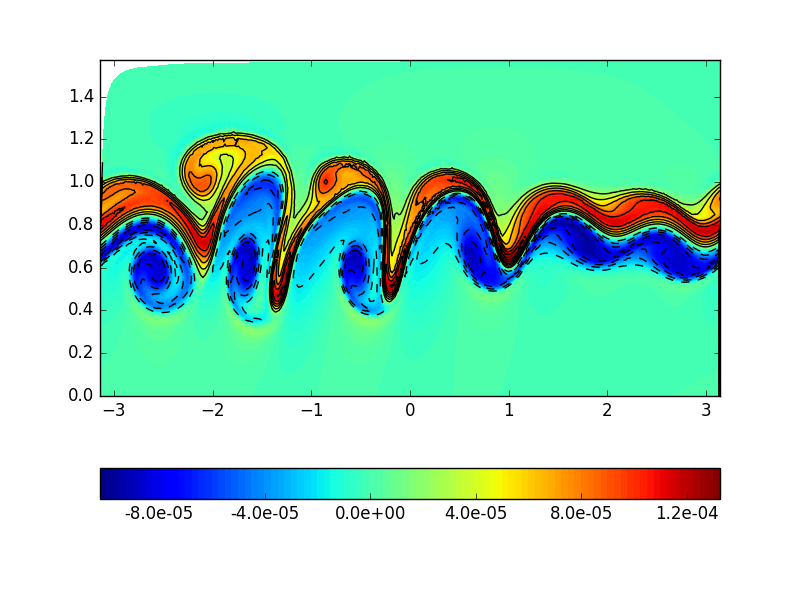}
\caption{As for Fig. \ref{fig::galewsky_1}, day 6.}
\label{fig::galewsky_3}
\end{center}
\end{figure}

The final test is for the nonlinear barotropic instability of a mid latitude jet \cite{GSP04}. The flow
is initialized with a near balanced state, perturbed by a shallow Gaussian hill in the depth field. As
the gravity wave triggered by this perturbation radiates outward it interacts with the mean flow, 
exciting a shear instability. At low resolutions, errors due to grid imprinting are of greater amplitude
than the Gaussian hill perturbation and a wave number 4 instability arises. At higher resolutions the 
Gaussian perturbation dominates the grid imprinting errors and an instability that matches well in position
and shape with the published results from a high resolution spectral model \cite{GSP04} emerges.
Results are presented for the vorticity field in figs. \ref{fig::galewsky_1}, \ref{fig::galewsky_2} and
 \ref{fig::galewsky_3} at days 4, 5, and 6 respectively with $p=3$, $N_e = 32$ and $\Delta t = 40\mathrm{s}$. 
The biharmonic viscosity is applied in order to suppress instabilities due to nonlinear cascades to the grid 
scale.

\section{Conclusions}\label{sec::conc}

A model of the rotating shallow water equations on a cubed sphere using mixed mimetic spectral elements 
is presented.
\secondRev{ The model preserves the optimal convergence of errors for the vorticity, velocity and the 
fluid depth on the $W_h$, $U_h$ and $Q_h$ function spaces respectively on the non-affine, smoothly 
varying mesh of the cubed sphere, as demonstrated for the standard test cases presented here.}

The mimetic properties of the incidence matrices are preserved independent of the geometry, ensuring the 
conservation of mass, vorticity and energy. While mass conservation holds to machine precision due to the 
point wise satisfaction of the divergence theorem, vorticity and energy conservation hold in the weak form.
As such conservation errors for vorticity remain bounded for the duration of the simulations presented,
and energy conservation errors converge as the product of the temporal and spatial orders of the scheme in 
the absence of viscosity. These results are validated using standard test cases.

\firstRev{
One curious result is that the $L_2$ errors for the divergence converge at one degree lower than their anticipated
rate, despite the fact that the errors for the fluid depth, which is also defined on $Q_h$, converge at the
correct rate. The reason for this is not well understood by the authors and requires further investigation.
}

\firstRev{
Preliminary performance results using 24 processors across 2 nodes shows that approximately 80\% of the
compute time is spent on matrix-matrix multiplication for local element operators during the assembly of
the nonlinear terms, and a further 5\% on the interpolation of vector fields to global coordinates
via the Piola transform. These results suggest that significant performance gains could be realized by
replacing the matrix multiplications for diagonal matrices by single loops over quadrature points.
}

In future work we intend to explore the extension of this method to the three dimensional primitive 
equations, as well as alternative formulations of the governing equations and the consequences of 
potential enstrophy conservation and exact quadrature via an iso-parametric Jacobian transformation that
is not reliant on the evaluation of transcendental functions.

\section{Acknowledgements}\label{sec::ack}

David Lee would like to thank Dr. Mark Taylor for several helpful discussions on the formulation of the 
Jacobian mapping, and Prof. Hugh Blackburn for the generous use of machine time for the production of the 
results. This research was supported as part of the Launching an Exascale ACME Prototype (LEAP) project, 
funded by the US Department of Energy, Office of Science, Office of Biological and Environmental Research.
This research used resources provided by the Los Alamos National Laboratory Institutional Computing Program, 
which is supported by the U.S. Department of Energy National Nuclear Security Administration under Contract 
No. DE-AC52-06NA25396.

%\section*{References}

%\bibliographystyle{elsarticle-num}
%\bibliography{references}

\end{document}